\documentclass{tran-l}
\usepackage{amssymb, graphicx}

\newtheorem{theorem}{Theorem}[section]
\newtheorem{lem}[theorem]{Lemma}
\newtheorem{cor}[theorem]{Corollary}
\newtheorem{prop}[theorem]{Proposition}

\theoremstyle{definition}
\newtheorem{defi}[theorem]{Definition}

\theoremstyle{remark}
\newtheorem{rem}[theorem]{Remark}

\numberwithin{equation}{section}

\newcommand{\al}{\alpha}
\newcommand{\be}{\beta}
\newcommand{\de}{\delta}

\newcommand{\si}{\sigma}

\newcommand{\Ga}{\Gamma}
\newcommand{\La}{\Lambda}
\newcommand{\Om}{\Omega}
\newcommand{\Si}{\Sigma}
\newcommand{\om}{\omega}
\renewcommand{\th}{\theta}
\newcommand{\ze}{\zeta}

\newcommand{\n}{\mathbf n}
\newcommand{\m}{\mathbf m}

\newcommand{\bZ}{\mathbf Z}
\newcommand{\CC}{\mathcal C}
\newcommand{\LL}{\mathcal L}
\newcommand{\EE}{\mathcal E}

\renewcommand{\SS}{\mathcal S}
\newcommand{\C}{\mathbb C}
\renewcommand{\L}{\mathbf L}

\newcommand{\Z}{\mathbb Z}

\newcommand{\R}{\mathbb R}

\newcommand{\CP}{{\mathbb C}{\mathbb P}}
\newcommand{\del}{\partial}
\newcommand{\x}{\times}
\newcommand{\lra}{\longrightarrow}
\newcommand{\hra}{\hookrightarrow}
\DeclareMathOperator{\lk}{lk}
\DeclareMathOperator{\tb}{tb}
\DeclareMathOperator{\rot}{rot}
\DeclareMathOperator{\PD}{PD}

\begin{document}

\title{On symplectic fillings of lens spaces} 
\author{Paolo Lisca}
\address{Dipartimento di Matematica ``L. Tonelli''\\ 
Universit\`a di Pisa \\I-56127
Pisa, ITALY} 
\email{lisca@dm.unipi.it}
\begin{abstract} 
Let $\overline\xi_{\rm st}$ be the contact structure naturally induced
on the lens space $L(p,q)=S^3/\Z/p\Z$ by the standard contact
structure $\xi_{\rm st}$ on the three--sphere $S^3$. We obtain a
complete classification of the symplectic fillings of
$(L(p,q),\overline\xi_{\rm st})$ up to orientation--preserving
diffeomorphisms. In view of our results, we formulate a conjecture on
the diffeomorphism types of the smoothings of complex two--dimensional
cyclic quotient singularities.
\end{abstract}
\maketitle

\section{Introduction and statement of results}
\label{s:intro}

A (coorientable) \emph{contact three--manifold} is a pair $(Y,\xi)$,
where $Y$ is a closed three--manifold and $\xi\subset TY$ a
two--dimensional distribution given as the kernel of a one--form
$\al\in\Om^1(Y)$ such that $\al\wedge d\al$ is a volume form. The
orientation on $Y$ determined by $\al\wedge d\al$ only depends on the
distribution $\xi$. We shall always assume that the underlying
manifold of a contact three--manifold is oriented, and that the
orientation is the one induced by the contact structure.

A \emph{symplectic filling} of a closed contact three--manifold
$(Y,\xi)$ is pair $(X,\om)$ consisting of a smooth, compact, connected
four--manifold $X$ and a symplectic form $\om$ on $X$ such that, if
$X$ is oriented by $\om\wedge\om$ and $\del X$ is given the boundary
orientation, there exists an orientation--preserving diffeomorphism
$\varphi\colon Y\to\del X$ such that $\om|_{d\varphi(\xi)}\neq 0$ at
every point of $\del X$. For example, the unit four--ball
$B^4\subset\C^2$ endowed with the restriction of the standard K\"ahler
form on $\C^2$ is a symplectic filling of $(S^3,\xi_{\rm st})$, where
the~\emph{standard contact structure} $\xi_{\rm st}$ on $S^3$ is, by
definition, the distribution of complex lines tangent to
$S^3\subset\C^2$.

The first classification result for symplectic fillings is due to
Eliashberg~\cite{El3}, who proved that if $(X,\om)$ is a symplectic
filling of $(S^3,\xi_{\rm st})$, then $X$ is diffeomorphic to a blowup
of $B^4$. McDuff~\cite{McD1} extended Eliashberg's result to the lens
spaces $L(p,1)$ endowed with the contact structure $\overline\xi_{\rm
st}$ defined in the following paragraph. Ohta and Ono~\cite{OO}
determined the diffeomorphism types of the strong symplectic fillings
of links of simple elliptic singularities endowed with their natural
contact structures. In~\cite{L3, OO2} results on the intersection
forms of symplectic fillings of finite quotients of $S^3$ are proved.

The distribution $\xi_{\rm st}$ on $S^3$, being invariant under the
natural action of $U(2)$, is a fortiori invariant under the induced
action of the subgroup
\begin{equation}\label{e:subgroup}
G_{p,q}=
\{\left(\begin{smallmatrix}
\ze & 0\\
0 & \ze^q
\end{smallmatrix}\right)
\ |\ \ze^p=1\}\subset U(2),
\end{equation}
where $p, q\in\Z$. It follows that when $p> q\geq 1$ and $p, q$ are
coprime, $\xi_{\rm st}$ descends to a contact structure
$\overline\xi_{\rm st}$ on the lens space $L(p,q)=S^3/G_{p,q}$. 

Let $D_p$ denote the disk bundle over the two--sphere with Euler class
$e(D_p)=p$. The contact three--manifold $(L(p,1),{\overline\xi_{\rm
st}})$ admits a symplectic filling of the form $(D_{-p},\om)$ for
every $p>1$. On the other hand, $(L(4,1),{\overline\xi_{\rm st}})$
also admits a symplectic filling of the form
$\left(\CP^2\setminus\nu(C),\om_0\right)$, where $\nu(C)$ is a
strictly pseudo--concave neighborhood of a smooth conic
$C\subset\CP^2$ and $\om_0$ is the restriction of the standard
K\"ahler form on $\CP^2$.

McDuff proved~\cite{McD1} that if $(X,\om)$ is a symplectic filling of
$(L(p,1),{\overline\xi_{\rm st}})$, then $X$ is
orientation preserving diffeomorphic to a smooth blowup of:
\begin{enumerate}
\item[(a)]
$D_{-p}$ if $p\neq 4$,
\item[(b)]
$D_{-4}$ or $\CP^2\setminus\nu(C)$ if $p=4$.
\end{enumerate}

In this paper we obtain a complete classification up to diffeomorphism
for the symplectic fillings of the contact three--manifolds
$\left(L(p,q),{\overline\xi_{\rm st}}\right)$. In order to state our
result we need to introduce some notation.

Let $p$ and $q$ be coprime integers such that
\begin{equation}\label{e:assumptions}
p>q\geq 1,\qquad
\frac p{p-q} = b_1 - \cfrac{1}{b_2 -
       \cfrac{1}{\ddots -
        \cfrac{1}{b_k}
}},\quad
b_1,\ldots,b_k\geq 2.
\end{equation}
The integers $b_i$ are uniquely determined by the rational number
$\frac p{p-q}$. Without the assumption $b_i\geq 2$, this uniqueness
property fails. The standard symbol $[b_1,\ldots, b_k]$ will be used
throughout the paper to denote a contined fraction as the one
in~\eqref{e:assumptions}.

{\bf Definition.}
A $k$--tuple of non--negative integers $(n_1,\ldots, n_k)$ 
is \emph{ admissible} if each of the denominators appearing 
in the continued fraction $[n_1,\ldots, n_k]$ is positive.

\begin{figure}[ht]
  \setlength{\unitlength}{1mm}
  \begin{center}
    \begin{picture}(100,22)
      \put(10,18){\large $n_1$}
      \put(29,18){\large $n_2$}
      \put(65,18){\large $n_{k-1}$}
      \put(86,18){\large $n_k$}
      \includegraphics[width=10cm]{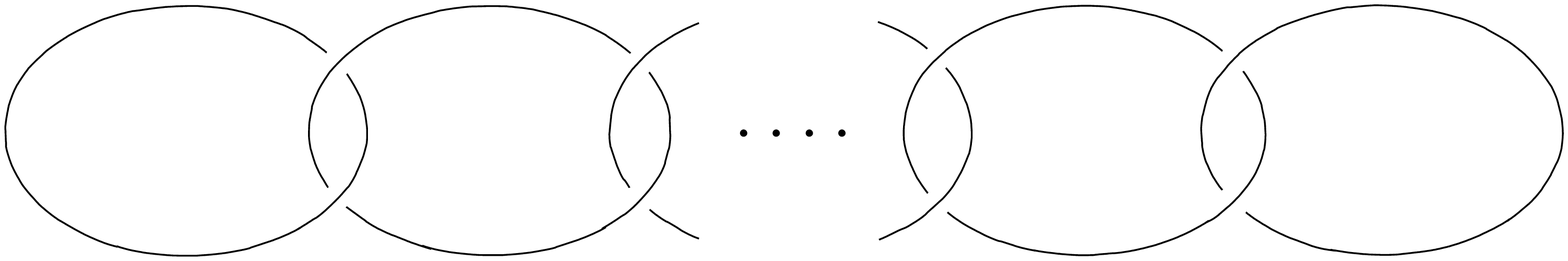}
    \end{picture}    
  \end{center}
\caption{The manifold $N(\n)$}
\label{f:L} 
\end{figure}

Let $\n=(n_1,\ldots, n_k)$ be an admissible $k$--tuple of
non--negative integers such that
\begin{equation}\label{e:vanishfrac}
[n_1,\ldots, n_k]=0.
\end{equation}
Let $N(\n)$ be the closed, oriented three--manifold given by surgery
on $S^3$ along the framed link of Figure~\ref{f:L}. It is easy to
check that Assumption~\eqref{e:vanishfrac} ensures the existence of an
orientation preserving diffeomorphism
\begin{equation}\label{e:varphi1}
\varphi\colon N(\n)\to S^1\x S^2.
\end{equation}

\begin{figure}[ht]
  \setlength{\unitlength}{1mm}
  \begin{center}
    \begin{picture}(100,35)
      \put(5.5,0){\large $b_1-n_1$}
      \put(25.5,0){\large $b_2-n_2$}
      \put(60.5,0){\large $b_{k-1}-n_{k-1}$}
      \put(83.5,0){\large $b_k-n_k$}
      \put(2,6){\small $-1$}
      \put(7.5,6){\small $-1$}
      \put(15,6){\small $-1$}
      \put(22,6){\small $-1$}
      \put(27.5,6){\small $-1$}
      \put(34.5,6){\small $-1$}
      \put(61,6){\small $-1$}
      \put(66,6){\small $-1$}
      \put(73,6){\small $-1$}
      \put(80,6){\small $-1$}
      \put(85,6){\small $-1$}
      \put(92,6){\small $-1$}
      \put(10,33){\large $n_1$}
      \put(29,33){\large $n_2$}
      \put(65,33){\large $n_{k-1}$}
      \put(86,33){\large $n_k$}
      \put(-2.5,11){$\L$}
      \put(12,12){\small ...}
      \put(31.5,12){\small ...}
      \put(70.5,12){\small ...}
      \put(89.5,12){\small ...}
      \put(0,3){\includegraphics[width=10cm]{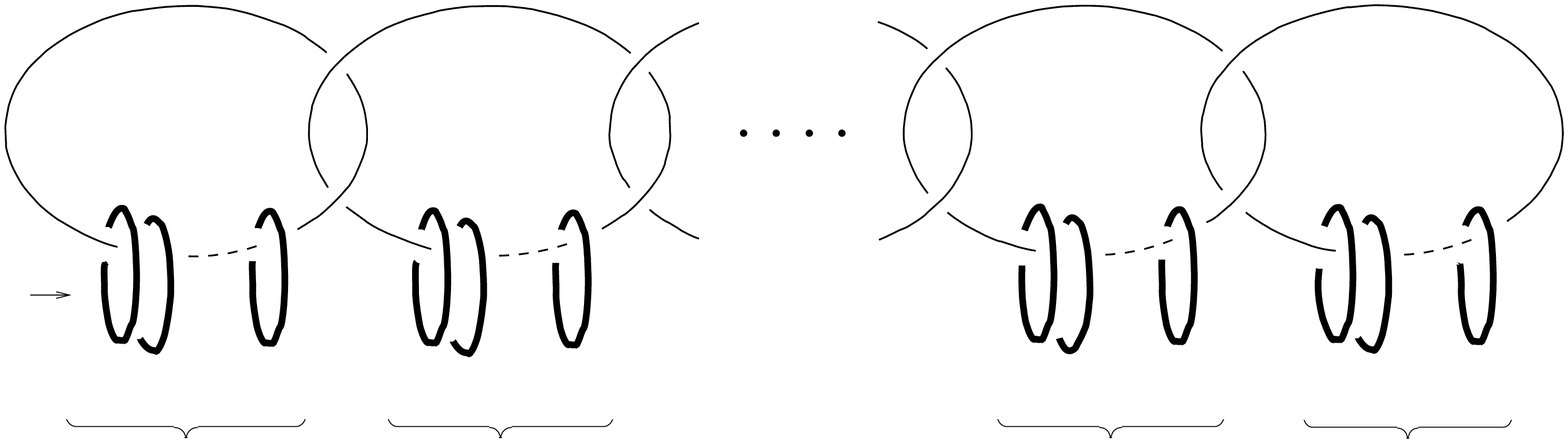}}
    \end{picture}    
  \end{center}
\caption{The framed link $\L\subset N(\n)$}
\label{f:W} 
\end{figure}

{\bf Definition.}
Fix a diffeomorphism $\varphi$ as in~\eqref{e:varphi1}, and let
${\mathbf L}\subset N(\n)$ be the thick framed link drawn in
Figure~\ref{f:W}. Define $W_{p,q}(\n)$ to be the smooth four--manifold
with boundary obtained by attaching two--handles to $S^1\x D^3$ along
the framed link
\[
\varphi({\mathbf L})\subset S^1\x S^2.
\]

{\bf Remark.}
The diffeomorphism type of $W_{p,q}(\n)$ is independent of the choice
of the diffeomorphism $\varphi$, because every self--diffeomorphism of
$S^1\x S^2$ extends to $S^1\x D^3$~\cite{Gl}.

{\bf Definition.} Let $\bZ_{p,q}\subset\Z^k$ be the set of
admissibible $k$--tuples of non--negative integers $(n_1,\ldots,n_k)$
such that
\[
[n_1,\ldots,n_k]=0\quad\text{and}\quad 
0\leq n_i\leq b_i\ \text{for}\ i=1,\ldots, k.
\]

It is easy to check (see Section~2) that the set $\bZ_{p,q}$ admits
the involution:
\[
\n=(n_1,\ldots,n_k)\mapsto\overline\n=(n_k,\ldots,n_1)
\]
Given coprime integers $p>q\geq 1$, we denote by $\overline
q$ the only integer satifsying
\[
p>{\overline q}\geq 1,\quad q{\overline q}\equiv 1\bmod p. 
\]
\begin{theorem}\label{t:main}
Let $p>q\geq 1$ be coprime integers. Then,
\begin{enumerate}
\item
Let $(W,\om)$ be a symplectic filling of the contact three--manifold
$(L(p,q),{\overline\xi_{\rm st}})$. Then, there exists
$\n\in\bZ_{p,q}$ such that $W$ is orientation preserving diffeomorphic
to a smooth blowup of $W_{p,q}(\n)$.
\item
For every $\n\in\bZ_{p,q}$, the four--manifold $W_{p,q}(\n)$ carries a
symplectic form $\om$ such that $\left(W_{p,q}(\n),\om\right)$ is a
symplectic filling of the contact three--manifold $\left(L(p,q),
{\overline\xi_{\rm st}}\right)$. Moreover, there are no classes in
$H_2(W_{p,q}(\n);\Z)$ with self--intersection equal to $-1$.
\item
Let $\n\in\bZ_{p,q}$ and $\n'\in\bZ_{p',q'}$. Then,
$W_{p,q}(\n)\# r\overline\CP^2$ is orientation preserving
diffeomorphic to $W_{p',q'}(\n')\# s\overline\CP^2$ if and only if: 
\begin{enumerate}
\item
$(p',s)=(p,r)$ and $(q',\n')=(q,\n)$, or 
\item
$(p',s)=(p,r)$ and $(q',\n')=({\overline q},{\overline\n})$.
\end{enumerate}
\end{enumerate}
\end{theorem}

Theorem~\ref{t:main} gives a complete diffeomorphism classification of
the symplectic fillings of $(L(p,q),\overline\xi_{\rm st})$, extending
the result of McDuff quoted above~\footnote{Portions (1) and (2) of
Theorem~\ref{t:main} were announced in~\cite{L1, L2}.}.  For instance,
as explained in~\cite{L1}, by Theorem~\ref{t:main} there are exactly
two symplectic fillings of $(L(p^2,p-1),{\overline\xi_{\rm st}})$ up
to blowups and diffeomorphisms. One comes from the canonical
resolution of the associated singularity, while the other is the
rational homology ball used in the symplectic rational blowdown
construction~\cite{FS, Sy}.

The following corollary gives three new results as applications of
Theorem~\ref{t:main}:

\begin{cor}\label{c:main}
Let $p>q\geq 1$ be coprime integers. 
\begin{enumerate}
\item[(a)] Given a positive integer $n$, there exist infinitely many
lens spaces $L(p,q)$ such that the contact three--manifold
$(L(p,q),\overline\xi_{\rm st})$ has more than $n$ symplectic fillings
pairwise distinct up to homotopy equivalence and whose underlying
four--manifolds are smoothly minimal~\footnote{A smooth, oriented
four--manifold $X$ is \emph{smoothly minimal} if the interior of $X$
contains no smoothly embedded two--sphere of self--intersection
$-1$.}.
\item[(b)] There exist infinitely many lens spaces $L(p,q)$ such that
$q>1$ and the contact three--manifold $(L(p,q),\overline\xi_{\rm st})$
has only one symplectic filling up to blowups and diffeomorphisms.
\item[(c)] The contact three--manifold $(L(p,q),\overline\xi_{\rm
st})$ has a symplectic filling $(W,\om)$ with $b_2(W)=0$ if and only
if $(p,q)=(m^2,mh-1)$, for some $m$ and $h$ coprime natural numbers.
\end{enumerate}
\end{cor}

The statement of Theorem~\ref{t:main} -- not the proof given here --
is related to the deformation theory of complex two--dimensional
cyclic quotient singularities. In fact, the contact three--manifold
$\left(L(p,q), {\overline\xi_{\rm st}}\right)$ can be viewed as the
link of the singularity $(\C^2/G_{p,q},0)$ together with the natural
contact structure given by the complex tangents. Every smoothing of
$(\C^2/G_{p,q},0)$ determines Stein fillings $F$ of
$(L(p,q),{\overline\xi}_0)$, and Theorem~\ref{t:main} implies that any
such $F$ must be diffeomorphic to $W_{p,q}(\n)$ for some
$\n\in\bZ_{p,q}$. It seems likely that a converse to this fact should
also hold, because every irreducible component of the reduced base
space $S_{\text{red}}$ of the versal deformation of $(\C^2/G_{p.q},0)$
gives a smoothing of the singularity, and Stevens~\cite{St} proved the
existence of a one--to--one correspondence between $\bZ_{p,q}$ and the
set ${\mathcal S}_{p,q}$ of irreducible components of $S_{\text{red}}$
(see also~\cite{Ch}). In light of the results obtained in this paper,
we propose the following:

{\bf Conjecture.} Let $F(\n)$ be a Stein filling of $\left(L(p,q),
{\overline\xi_{\rm st}}\right)$ determined by the smoothing of
$(\C^2/G_{p.q},0)$ corresponding to $\n\in\bZ_{p,q}$ under the Stevens
correspondence~\cite{St}. Then, $F(\n)$ is diffeomorphic to
$W_{p,q}(\n)$.

There is evidence supporting this conjecture. In fact, in this paper
 we prove that each $W_{p,q}(\n)$ carries Stein structures
 (Corollary~\ref{c:stein}). Moreover, the smoothing corresponding to
 $(1,2,\ldots,2,1)\in\bZ_{p,q}$ is known to be isomorphic to the
 canonical resolution $X_{p,q}$ of the singularity, and it is not hard
 to verify that if $\n=(1,2,\ldots, 2,1)$ then $W_{p,q}(\n)$ is indeed
 diffeomorphic to a regular neighborhood of the exceptional divisor in
 $X_{p,q}$.  Finally, by~\cite[5.9.1]{Wa}, a singularity
 $(\C^2/G_{p.q},0)$ has a smoothing with $b_2=0$ if and only if
 $(p,q)=(m^2,mh-1)$, with $m$ and $h$ coprime, in agreement with
 Corollary~\ref{c:main}(c).

The paper is organized as follows. In Section~\ref{s:corollary} we
prove Corollary~\ref{c:main} assuming Theorem~\ref{t:main}. In
Section~\ref{s:gluing} we show that every symplectic filling $(W,\om)$
of $(L(p,q),\overline\xi_{\rm st})$ can be compactified to a rational
symplectic four--manifold $X$ so that $X\setminus W$ is a neighborhood
of an immersed symplectic surface $\Ga\subset X$ of a special kind. In
Section~\ref{s:strings} we prove Theorem~\ref{t:main}(1). In
Section~\ref{s:stein} we construct Stein structures on the smooth
four--manifolds $W_{p,q}(\n)$. In Section~\ref{s:contact} we prove
Theorem~\ref{t:main}(2). In Section~\ref{s:mainproof} we complete the
proof of Theorem~\ref{t:main}.

{\bf Acknowledgements.}  The author gratefully acknowledges support
and hospitality from the Department of Mathematics of the University
of Georgia during the preparation of part of this paper. Warm thanks
go to Patrick Popescu--Pampu for his interest in my work and for
several remarks and suggestions which allowed me to improve this paper
in a few points where it lacked precision, and to Andr\'as Stipsicz
for useful comments.

\section{The proof of Corollary~\ref{c:main}}
\label{s:corollary}

\begin{defi}\label{d:blowup}
Given a $k$--tuple of positive integers
\[
(n_1,\ldots, n_{s-1}, n_s, n_{s-1},\ldots, n_k)
\]
with $n_s=1$, we say that the $(k-1)$--tuple
\[
(n_1,\ldots, n_{s-1}-1, n_{s+1}-1,\ldots, n_k)
\]
is obtained by a {\em blowdown at $n_s$} (with the obvious meaning
of the notation when $s=1$ or $s=k$). The reverse process is a {\em
blowup}. A blowdown at $n_s$ is {\em strict} if $s>1$. The reverse
process is a {\em strict} blowup.
\end{defi}

It is showed in~\cite[Appendix]{OW} that a $k$--tuple of positive
integers $(n_1,\ldots,n_k)$ is admissible~\footnote{The definition of
admissibility given in~\cite[Appendix]{OW} is easily seen to be
equivalent to the one given in Section~1.} if and only if the symmetric
matrix
\[ 
\begin{pmatrix}
n_1 & -1 & & &\\
-1 & n_2 & -1 & &\\
  & -1 & n_3 & &\\
& & & \ \ddots & -1\\
& & & -1 &-n_k\\
\end{pmatrix}
\]
is positive semi--definite of rank $\geq k-1$. It immediately follows
from this fact that if $(n_1,\ldots,n_k)$ is admissible then
$(n_k,\ldots,n_1)$ is admissible, 
\[
(n_i,n_{i+1},\ldots,n_{j-1},n_j)
\]
is admissible for every $1\leq i\leq j\leq k$, and blowing up and
blowing down preserve admissibility.

\begin{lem}\label{l:0-sequence}
Let $(n_1,\ldots, n_k)$ be a $k$--tuple of positive integers. Then, 
the following conditions are equivalent:
\begin{enumerate}  
\item[$\bullet$]
$(n_1,\ldots, n_k)$ is admissible and satisfies $[n_1,\ldots, n_k]=0$; 
\item[$\bullet$]
$(n_1,\ldots, n_k)$ is obtained from $(0)$ by a sequence of strict blowups.
\end{enumerate}
\end{lem}

\begin{proof}
Clearly, a $k$--tuple of positive integers is obtained from $(0)$ via
a sequence of strict blowups if and only if it is obtained from
$(1,1)$ via such a sequence. Moreover, both the property of being
admissible and that of having vanishing associated continued fraction
are preserved under blowup. Therefore, since $(1,1)$ is admissible and
$[1,1]=0$, any $k$--tuple of positive integers obtained from $(0)$ by
a sequence of strict blowups is admissible and has vanishing
associated continued fraction.

Conversely, let $(n_1,\ldots,n_k)$ be an admissible $k$--tuple of
positive integers with $[n_1,\ldots,n_k]=0$. Then, we must have $k\geq
2$. We shall argue by induction that $(n_1,\ldots,n_k)$ is obtained
from $(0)$ by a sequence of strict blowups. For $k=2$ the statement is
obvious, so suppose that $k>2$ and the statement corresponding to
$k-1$ holds true. An easy induction argument shows that if $n_i\geq 2$
for every $i=1,\ldots, k$, then $[n_1,\ldots, n_k] > 1$. Thus, if
$[n_1,\ldots, n_k]=0$ then necessarily $n_i=1$ for some
$i\in\{1,\ldots, k\}$. We conclude that $(n_1,\ldots,n_k)$ is obtained
from the admissible $(k-1)$--tuple
\[
(n_1,\ldots,n_{i-1}-1,n_{i+1}-1,\ldots,n_k)
\]
by a blowup (which is strict if and only if $i>1$). By the induction
hypothesis,
\[
(n_1,\ldots,n_{i-1}-1,n_{i+1}-1,\ldots,n_k)
\]
is obtained from $(0)$ by a sequence of strict blowups, hence 
$n_j=1$ for some $j>1$. Therefore, the $k$--tuple
$(n_1,\ldots,n_k)$ is obtained by a strict blowup from the
admissible $(k-1)$--tuple 
\[
(n_1,\ldots,n_{j-1}-1,n_{j+1}-1,\ldots,n_k), 
\]
and the conclusion follows by induction.
\end{proof}

\begin{lem}\label{l:uniqueness}
Let $p>q\geq 1$ be coprime integers, and suppose that 
\[
\frac pq=[a_1,\ldots,a_h],\quad\text{with}\quad a_1,\ldots,a_h\geq 5.
\]
Then, 
\[
\bZ_{p,q}=\{(1,2,\ldots,2,1)\}.
\]
\end{lem}

\begin{proof}
Let 
\[
\frac p{p-q}=[b_1,\ldots,b_k],\quad b_1,\ldots, b_k\geq 2.
\]
Using Riemenschneider's point diagram~\cite{Ri}, one can easily check
that the assumption on the $a_i$'s implies the following three
conditions:
\begin{itemize}
\item
$k\geq 4$, 
\item
$b_1,\ldots, b_k\leq 3$,
\item
if $b_i=b_j=3$ then either $3<i=j<k-2$ or $|i-j|\geq 3$.
\end{itemize}
Therefore, if $(n_1,\ldots,n_k)\in\bZ_{p,q}$ we have
\begin{enumerate}
\item
$k\geq 4$,
\item
$n_i\leq 3$ for every $i=1,\ldots,k$,
\item
if $n_i=n_j=3$, then either $3<i=j<k-2$ or $|i-j|\geq 3$.
\end{enumerate}

We shall argue by induction on $k\geq 4$ that if $(n_1,\ldots,n_k)$ is
an admissible sequence of non--negative integers such that
$[n_1,\ldots,n_k]=0$ and such that (1), (2) and (3) above hold, then
\[
(n_1,\ldots,n_k)=(1,2,\ldots,2,1).
\]

If $k=4$ one
immediately sees that, by Lemma~\ref{l:0-sequence} and Assumption (3), 
\[
(n_1,n_2,n_3,n_4)=(1,2,2,1).
\]
Now suppose $k>4$. By Lemma~\ref{l:0-sequence}, $n_i>0$ for every
$i=1,\ldots,k$ and $n_j=1$ for some index $j>1$. We claim that
$j=k$. In fact, suppose $j<k$. If $j=2$ or $j=k-1$, then Assumption
(3) implies $n_{j-1}=n_{j+1}=2$. By Lemma~\ref{l:0-sequence} this is
impossible, because two blowdowns would give a string of
length bigger than $1$ containing $0$. Therefore we have
$2<j<k-1$. Blowing down once yields the sequence
\[
(n_1,\ldots,n_{j-1}-1,n_{j+1}-1,\ldots, n_k),
\]
which still satisfies the three assumptions. Therefore, by 
induction
\[
(n_1,\ldots,n_{j-1}-1,n_{j+1}-1,\ldots, n_k)=(1,2,\ldots,2,1)
\]
But then $n_{j-1}=n_{j+1}=3$, which goes against Assumption (3). We
conclude that $j=k$. Blowing down once yields the sequence
\[
(n_1,\ldots,n_{k-1}-1), 
\]
which satisfies the assumptions if and only if $n_{k-3}<3$. Notice
that if $n_{k-3}<3$, we can apply induction and reach the
conclusion. Blowing down two more times yields the sequence
\[
(n_1,\ldots, n_{k-3}-1),
\]
which satisfies the assumptions, so by induction we get $n_{k-3}=2$.
\end{proof}

We need the following elementary facts about continued fractions
(see~e.g.~\cite[Appendix]{OW} for the proofs):
\begin{itemize}
\item
Let $p>q\geq 1$ be coprime integers, and suppose that 
\[
\frac pq=[a_1,a_2,\ldots,a_h],\quad a_1,\ldots,a_h\geq 2.
\]
Then, 
\[
[a_h,a_{h-1},\ldots,a_1]=\frac p{\overline q},
\]
where $p>{\overline q}\geq 1$ and $q{\overline q}\equiv 1\bmod p$. 
\item
If $(n_1,\ldots, n_k)$ is an admissible $k$--tuple of positive integers, then
\[
[n_1,n_2,\ldots,n_k]=0
\]
if and only if 
\[
[n_k,n_{k-1}\ldots,n_1]=0.
\] 
\end{itemize}

\begin{lem}\label{l:0-frac-induction}
Let $(n_1,\ldots,n_k)$ be an admissible $k$--tuple of
positive integers such that
\[
[n_1,n_2,\ldots,n_k]=0,\quad k\geq 3.
\]
Suppose that there is exactly one index $j\in\{1,\ldots,k\}$ such that
$n_j=1$.
Then, 
there are coprime integers $m$ and $n$ such that
\[
[n_1,n_2,\ldots,n_{j-1},2,n_{j+1},\ldots,n_k]
=\frac{m^2}{mn+1}.
\]
\end{lem}

\begin{proof}
We argue by induction on $k\geq 3$. For $k=3$, the assumptions and
Lemma~\ref{l:0-sequence} imply
\[
(n_1,n_2,n_3)=(2,1,2),
\]
and 
\[
[2,2,2]=\frac 43
\]
is of the stated form. Now suppose $k>3$. Since, by
Lemma~\ref{l:0-sequence}, $(n_1,\ldots,n_k)$ must be obtained from
$(2,1,2)$ by a sequence of strict blowups, we have 
\[
n_1=2,\  n_k>2,\quad\text{or}\quad 
n_1>2,\  n_k=2.
\]
Observe that
\[
[n_1,n_2,\ldots,n_{j-1},2,n_{j+1},\ldots,n_k]
=\frac{m^2}{mn+1}
\]
if and only if 
\[
[n_k,n_{k-1},\ldots,n_{j+1},2,n_{j-1},\ldots,n_1]
=\frac{m^2}{m(m-n)+1}.
\]
Therefore, without loss of generality we may assume 
\[
n_1=2\quad\text{and}\quad n_k>2. 
\]
Under these assumptions we claim that
\begin{equation}\label{e:claim}
[n_2,\ldots,n_k-1]=0.
\end{equation}
In order to prove the claim, we argue by induction on $k\geq 4$. If $k=4$
we must have
\[
(n_1,n_2,n_3,n_4)=(2,2,1,3),
\]
and the claim is clear. If $k>4$, observe that, since by
Lemma~\ref{l:0-sequence}, $(n_1,\ldots, n_k)$ is a blowup of
$(2,2,1,3)$, we must have $j>2$. Then, blowing down once yields the
string
\[
(n_1=2,\ldots,n_{j-1}-1,n_{j+1}-1,\ldots,n_k),
\]
to which we may apply induction to conclude
\[
[n_2,\ldots,n_{j-1}-1,n_{j+1}-1,\ldots,n_k-1]=0 
\]
if $j<k-1$, and 
\[
[n_2,\ldots,n_{k-2}-1,n_k-2]=0
\]
if $j=k-1$. Blowing up again proves the claim~\eqref{e:claim}. 

Now induction applied to~\eqref{e:claim} gives 
\[
[n_2,\ldots,n_{j-1},2,n_{j+1},\ldots,n_k-1]=
\frac{n^2}{nh+1},
\]
with $n$, $h$ coprime integers. Since
\[
(1+nh)(1-nh)\equiv 1\bmod n^2,
\]
we have
\[
[n_k-1,\ldots,n_{j+1},2,n_{j-1},\ldots,n_2]=
\frac{n^2}{n(n-h)+1},
\]
therefore
\[
[n_k,\ldots,n_{j+1},2,n_{j-1},\ldots,n_2]=
\frac{n^2}{n(n-h)+1} +1 = \frac{2n^2 -nh+1}{n(n-h)+1}.
\]
Thus, since 
\[
(n^2-nh+1)(2nh-h^2+2)=(2n^2-nh+1)(nh-h^2+1)+1,
\]
we conclude 
\[
[n_1=2,n_2,\ldots,n_{j-1},2,n_{j+1},\ldots,n_k]=
2-\frac{2(nh+1)-h^2}{2n^2-nh+1}=
\frac{m^2}{mn+1},
\]
where $m=2n-h$.
\end{proof}

\begin{proof}[Proof of Corollary~\ref{c:main}]
(a) Let $p>q\geq 1$ be coprime and such that 
\[
\frac p{p-q}=[b_1,\ldots,b_k],\quad b_1,\ldots,b_k\geq 2,
\]
with $k\geq 4$, $b_2,\ldots,b_{k-2}\geq 3$ and $b_k\geq k-2$. Let
$r\in\Z$, $0\leq r\leq k-4$. Then, 
\[
\n_r=(1,\overbrace{2,\ldots,2}^r,3,\overbrace{2,\ldots,2}^{k-4-r},1,k-2-r)
\in \bZ_{p,q}
\]
One can easily check that
\begin{equation}\label{e:rank}
\chi\left(W_{p,q}(\n_r)\right)
 = 5+ \sum_{i=1}^k (b_i - 3) + r.
\end{equation}
Then, Equation~\eqref{e:rank} and Theorem~\ref{t:main} imply that
$(L(p,q),{\overline\xi_{\rm st}})$ admits at least $k-3$ smoothly minimal
symplectic fillings up to homotopy equivalence.

(b) If 
\[
\frac pq=[a_1,\ldots,a_h],\quad a_1,\ldots,a_h\geq 5,
\]
then by Lemma~\ref{l:uniqueness} 
\[
\bZ_{p,q}=\{(1,2,\ldots,2,1)\}.
\]
The conclusion follows by Theorem~\ref{t:main}.

(c) Suppose that
\[
\frac {p}{p-q}=[b_1,\ldots,b_k].
\]
It is easy to check the attaching circle of each two--handle
of $W_{p,q}(\n)$ is homologically non--trivial in $S^1\x
S^2$. Therefore, $b_2(W_{p,q}(\n))=0$ if and only if there is exactly
one index $j\in\{1,\ldots,k\}$ such that $n_j=1$, and
\[
(b_1,\ldots,b_k)=(n_1,n_2,\ldots,n_{j-1},2,n_{j+1},\ldots,n_k).
\]
Then, by Lemma~\ref{l:0-frac-induction},
\[
\frac{p}{p-q}=\frac{m^2}{mn+1}, 
\]
with $m$ and $n$ coprime. Therefore $q=mh-1$, with $h=m-n$.
\end{proof}

\section{Compactifications of symplectic fillings}
\label{s:gluing}

The purpose of this section is to prove Theorem~\ref{t:string} below.
In order to state the theorem, we need a definition.

\begin{defi}\label{d:s-string}
Let $(X,\om)$ be a symplectic four--manifold. A \emph{symplectic
string} in $X$ is an immersed symplectic surface
\[
\Ga=C_0\cup C_1\cup\cdots\cup C_k\subset X,
\]
where:
\begin{itemize}
\item
$C_i$ is a connected, embedded symplectic sphere for $i=0,\ldots, k$;
\item 
$C_i$ and $C_{i+1}$ intersect transversely and positively at a single point, 
for $i=0,\ldots, k-1$;
\item
$C_i\cap C_j=\emptyset$ if $|i-j|>1$, for $i,j=0,\ldots,k$.
\end{itemize}
We say that $\Ga$ as above is \emph{of type $(m_0,\ldots,m_k)$}
if, furthermore,
\begin{itemize}
\item
$C_i\cdot C_i = m_i$ for $i=0,\ldots, k$.
\end{itemize} 
\end{defi}

\begin{theorem}\label{t:string}
Let $p>q\geq 1$ be coprime integers, and suppose that 
\[
\frac p{p-q}=[b_1,\ldots,b_k],\quad b_1,\ldots,b_k\geq 2.
\]
Let $(W,\om)$ be a symplectic filling of $(L(p,q),{\overline\xi_{\rm
st}})$.  Then, for some integer $M\geq 0$, $W$ is
orientation preserving diffeomorphic to the complement of a regular
neighborhood of a symplectic string
\[
\Ga=C_0\cup C_1\cup\cdots\cup C_k\subset \CP^2\# M{\overline \CP}^2
\]
of type $(1,1-b_1,-b_2,\ldots,-b_k)$, where $\CP^2\# M{\overline
\CP}^2$ is endowed with a symplectic blowup $\om_M$ of the standard
K\"ahler form on $\CP^2$, and $C_0$ is a complex line in $\CP^2$.
\end{theorem}

Theorem~\ref{t:string} will be used in Section~\ref{s:strings}. We
start with an auxiliary construction of a suitable symplectic form on
$(0,\infty)\x L(p,q)$.

The function
\begin{eqnarray*}
\rho\colon & \C^2 & \longrightarrow (0,+\infty)\\ 
& (z_1,z_2) & \mapsto |z_1|^2 + |z_2|^2
\end{eqnarray*}
is a K\"ahler potential for the standard symplectic form $\om_0$, i.e.
\[
\om_0=\frac i2\sum_{k=1}^2 dz_k\wedge d{\overline z}_k
=\frac i2 \del{\overline\del}\rho.
\]
Let $J_0$ be the standard complex structure on $\C^2$, and 
consider the one--form 
\begin{equation*}
\si_0 = -\frac 14 d\rho\circ J_0=-\frac i4 (\del\rho-\overline\del\rho)
=-\frac i4 \sum_{k=1}^2 (\overline z_k dz_k - z_k d\overline z_k).
\end{equation*}
Let $i\colon S^3=\rho^{-1}(1)\hookrightarrow\C^2$ be the inclusion map,
and define $\al_0=i^*\si_0\in\Om^1(S^3)$. Since the standard contact
structure $\xi_{\rm st}$ on $S^3$ is given by complex tangent lines,
we have $\xi_{\rm st}=\{\al_0=0\}\subset TS^3$.  Define
$\pi\colon\C^2\setminus\{(0,0)\} \to S^3$ by setting:
\[
\pi(\mathbf z)= 
\frac {\mathbf z}{\rho(\mathbf z)^{\frac 12}}. 
\]
The pair $(\rho,\pi)$ induces an orientation preserving
diffeomorphism:
\begin{equation*}
(\rho,\pi)\colon\C^2\setminus\{(0,0)\} \to (0,+\infty)\x S^3.
\end{equation*}
The diffeomorphism $(\rho,\pi)$ sends the standard symplectic form
$\om_0$ to the form $d(t\al_0)$, where $t$ is the coordinate on the
first factor. To see this, notice that
\[
(i\circ\pi)^*\si_0 = -\frac i4 \sum_{k=1}^2 \left(\overline z_k d
\left(\frac {z_k}{\rho(\mathbf z)^{\frac 12}}\right) - z_k
d\left(\frac {\overline z_k}{\rho(\mathbf z)^{\frac 12}}\right)\right)
= \frac 1\rho \si_0,
\]
and therefore:
\begin{equation}\label{e:omega0}
(\rho,\pi)^* d(t\al_0) = d(\rho \pi^* i^*\si_0) = d\si_0 =
-\frac i4 \left({\overline\del}\del\rho - \del{\overline
\del}\rho\right) = \frac i2\del{\overline\del}\rho.
\end{equation}
Since $U(2)$ acts on $\C^2$ via norm--preserving complex linear
transformations, $\si_0$ and $\om_0$ are clearly $U(2)$--invariant
while $\al_0$ is invariant under the naturally induced action on
$S^3$. Moreover, $(\rho,\pi)$ is $U(2)$--equivariant in an obvious
sense, so we have: 
\begin{itemize}
\item
a symplectic form ${\overline\om_0}$ induced by $\om_0$ on
$\C^2\setminus\{(0,0)\}/G_{p,q}$, where $G_{p,q}$ is the subgroup of
$U(2)$ defined by~\eqref{e:subgroup},
\item
a contact form $\overline\al_0$ induced by $\al_0$ on
$L(p,q)=S^3/G_{p,q}$, and
\item
a symplectomorphism
\begin{equation}\label{e:symplect}
\left(\C^2\setminus\{(0,0)\}/G_{p,q}, {\overline\om_0}\right)\cong
\left((0,+\infty)\x L(p,q), d(t{\overline\al_0})\right).
\end{equation}
\end{itemize}
With the notation of Section~\ref{s:intro}, we have
${\overline\xi_{\rm st}}=\{\overline\al_0=0\}$. The action
\begin{eqnarray}\label{e:S3-action}
S^1\x S^3 & \lra & S^3\\
(e^{i\th},{\mathbf z}) & \longmapsto & e^{i\th}{\mathbf z}
\end{eqnarray}
commutes with the action of $G_{p,q}$, so it induces a fixed--point
free $S^1$--action: 
\begin{equation}\label{e:standardaction}
S^1\x L(p,q) \lra L(p,q).
\end{equation}
There is also an obviously induced $S^1$--action on $(0,+\infty)\x
L(p,q)$. Since $\al_0$ is $U(2)$--invariant, the symplectic form
$d(t{\overline\al_0})$ is invariant under this $S^1$--action.

The hypersurfaces 
\[
\{t\}\x L(p,q)\subset (0,+\infty)\x L(p,q)
\]
satisfy an important property with respect to the $S^1$--invariant
symplectic form $d(t{\overline\al_0})$. We shall establish this
property after recalling its general context.

Let $(X,\om)$ be a symplectic four--manifold and $Y\subset X$ a
separating hypersurface endowed with a fixed--point free
$S^1$--action.  Following~\cite{McW}, we say that $Y$ is {\em
$\om$--compatible} if the orbits of the $S^1$--action lie in the null
directions of $\om|_Y$. In general, we say that an embedding $j\colon
Y\hra X$ is $\om$--compatible if $j(Y)\subset X$ is a separating
$\om$--compatible hypersurface.

An $\om$--compatible hypersurface $Y$ has a canonical
co--orientation. In fact, if $N\in T_p Y$, $p\in Y$, is the vector
field generating the $S^1$--action on $Y$ then, for every vector $V\in
T_p X$ transverse to $Y$, $\om(N,V)\neq 0$, for otherwise $\om$ would
be degenerate on $T_p X$. Thus,
\[
T_p X\setminus T_p Y = T_p^+ X \cup T_p^- X,
\]
where 
\[
T_p^+ X=\{V\in T_p X\ |\ \om(N,V)>0\},\quad 
T_p^- X=\{V\in T_pX \ |\ \om(N,V)<0\}. 
\]
This allows one to distinguish the two components of $X\setminus Y$ by
setting $X^+$ to be the component with inward pointing normal vector
$V\in T_p^+ X$, for all $p\in Y$, and $X^-$ to be remaining
component. Thus, we have:
\begin{equation*}
X=X^-\cup_Y X^+.
\end{equation*}

The importance of the $\om$--compatibility condition is due to the
fact that symplectic four--manifolds can be glued along
$\om$--compatible hypersurfaces. More precisely, J.~McCarthy and
J.~Wolfson proved the following result:

\begin{theorem}[\cite{McW}, Theorem~4.2]\label{t:mcw}
Let $Y$ be closed, oriented three--manifold with a fixed--point free
$S^1$--action. Let $(X_i,\om_i)$, $i=1,2$, be symplectic
four--manifolds and $j_i\colon Y\to X_i$, $i=1,2$ $\om$--compatible
embeddings. Then, there is a symplectic structure $\om$ on a smooth
four--manifold
\begin{equation*}
X=X_1^-\cup_Y X_2^+
\end{equation*}
obtained by gluing $X_1^-$ to $X_2^+$ along $Y$. Moreover, there are
neighborhoods $\nu_i(j_i(Y))$ of $Y$ in $X_i$, $i=1, 2$, such that
\[
\om|_{X^+_2\setminus \nu_2(Y)}=\om_2\quad\text{and}\quad 
\om|_{X_1^-\setminus \nu_1(Y)}=c\om_1
\]
for some constant $c>0$.\qed
\end{theorem}

The idea of the proof of Theorem~\ref{t:string} is to show that any
symplectic filling $(W,\om)$ of $(L(p,q),\overline\xi_{\rm st})$ can
be compactified in a suitable way. In order to do that, we first need
to establish the following:

\begin{prop}\label{p:gluing}
Let $(W, \om)$ be a symplectic filling of $(L(p,q),{\overline\xi_{\rm
st}})$.  Let $(X,\eta)$ be a symplectic four--manifold and $j\colon
L(p,q)\to X$ an embedding which is $\eta$--compatible with respect to
the $S^1$--action~\eqref{e:standardaction}. Then, there is a symplectic
form $\overline\om$ on a smooth four--manifold
\begin{equation*}
Z=W\cup_{\del W=\del X^-} X^-
\end{equation*}
obtained by gluing $X^-$ to $W$ along their boundaries.
\end{prop}

Proposition~\ref{p:gluing} will be proved after the next two lemmas.

\begin{lem}\label{l:compatible}
For every $t\in (0,+\infty)$, the hypersurface
\begin{equation*}
\{t\}\x L(p,q)\subset (0,+\infty)\x L(p,q)\quad 
\end{equation*}
is $S^1$--invariant and $d(t{\overline\al_0})$--compatible with
respect to the action~\eqref{e:standardaction}. Moreover,
\begin{align*}
\left[(0,+\infty)\x L(p,q)\right]^+ =&
(0,t]\x L(p,q),\\
\left[(0,+\infty)\x L(p,q)\right]^- =&
[t,+\infty)\x L(p,q)
\end{align*}
\end{lem}

\begin{proof}
The hypersurface $\{t\}\x L(p,q)$ is clearly $S^1$--invariant for
every $t\in (0,+\infty)$. In order to see that it is also
$d(t{\overline\al_0})$--compatible, we first show that the
$S^1$--invariant hypersurface 
\[
\rho^{-1}(t)\subset \C^2\setminus\{(0,0)\} 
\]
is $\om_0$--compatible. Given ${\mathbf z}\in\rho^{-1}(t)$, we have
\begin{equation*}
{\frac d{d\th}}{\Big |}_{\th=0}\left(e^{i\th} {\mathbf z}\right)
= i {\mathbf z}
\end{equation*}
Equivalently, the vector field
\begin{equation}\label{e:vfield}
N = J_0 \left(z_1\frac\del{\del z_1} + z_2\frac\del{\del z_2}\right)
\end{equation}
generates the $S^1$--action. Since $\si_0$ is
invariant,
\begin{equation}\label{e:action}
0 = L_N \si_0 = d i_N\si_0 + i_N d\si_0,
\end{equation}
and since $d\rho =\sum_k {\overline z_k}dz_k + z_k d{\overline z_k}$, we
have
\begin{equation*}
i_N \si_0  = -\frac 14 d\rho (J_0(N)) 
= \frac\rho 4
\end{equation*}
Therefore, by~\eqref{e:action} we have 
\begin{equation}\label{e:kernel}
i_N d\si_0 |_{\rho^{-1}(t)} = i_N \om_0 |_{\rho^{-1}(t)} = 0.
\end{equation}
Equation~\eqref{e:kernel} says that
\[
\rho^{-1}(t)\subset\C^2\setminus\{(0,0)\}
\]
is $\om_0$--compatible. Therefore, by $U(2)$--invariance
and~\eqref{e:symplect} the hypersurface 
\[
\{t\}\x L(p,q)\subset (0,+\infty)\x L(p,q)
\]
is $d(t{\overline\al_0})$--compatible.

In order to prove the second part of the statement, by
$U(2)$--invariance it is enough to check that, with respect to the
$S^1$--action~\eqref{e:S3-action} and the $d(t\al_0)$--compatible
hypersurface 
\[
\{t\}\x S^3\subset (0,+\infty)\x S^3, 
\]
we have
\begin{equation*}
\left[(0,+\infty)\x S^3\right]^+ = (0,t]\x S^3,\qquad
\left[(0,+\infty)\x S^3\right]^- = [t,+\infty)\x S^3.
\end{equation*}
In fact, the vector field 
\[
V = z_1\frac\del{\del z_1} + z_2\frac\del{\del z_2}
\]
is transverse to every hypersurface $\{t\}\x S^3$ and points inwardly
with respect to 
\[
[t,+\infty)\x S^3.
\]
Since $N=J_0 V$, we have
\begin{equation*}
\om_0(N, V) = \om_0(J_0 V, V) = -\om_0(V, J_0 V) < 0.
\end{equation*}
This concludes the proof.
\end{proof}

\begin{lem}\label{l:end}
Let $(W,\om)$ be a symplectic filling of $(L(p,q),{\overline\xi_{\rm
st}})$.  Then, there exists a symplectic four--manifold $(\widetilde
W,\widetilde\om)$ with one end $\EE$ such that:
\begin{enumerate}
\item[(a)] $(\EE,\widetilde\om|_{\EE})$ is symplectomorphic to
\[
\left((D,+\infty)\x L(p,q), d(t{\overline \al_0})\right)
\]
for some constant $D>0$;
\item[(b)] ${\widetilde W}\setminus\EE$ is orientation preserving
diffeomorphic to $W$.
\end{enumerate}
\end{lem}

\begin{proof}
Choose an open collar $C\subset W$ around $\del W$, and let 
\[
\varphi\colon C\setminus\del W\cong (0,+\infty)\x L(p,q)
\]
be an orientation preserving diffeomorphism such that, if we keep
denoting by $\om$ the push--forward of $\om$ to
\[
(0,+\infty)\x L(p,q)
\]
and we identify $\{1\}\x L(p,q)$ with $L(p,q)$, we have
$\om|_{\overline\xi_{\rm st}}\neq 0$. In fact, we claim that there
exists a contactomorphism
\[
\overline c\colon (L(p,q),\overline\xi_{\rm st})\to 
(L(p,q),\overline\xi_{\rm st})
\]
which reverses the orientation of $\overline\xi_{\rm st}$. Consequently, 
we can choose the identification of $\{1\}\x L(p,q)$ with $L(p,q)$ so that 
$\om|_{\overline\xi_{\rm st}}> 0$. 

In order to prove the claim, recall that  the standard contact structure 
$\xi_{\rm st}$ on $S^3$ has a natural orientation induced by the standard 
complex structure on $\C^2$. The orientation on $\xi_{\rm st}$ induces an 
orientation on $\overline\xi_{\rm st}$. Moreover, the map $\C^2\to\C^2$ 
given by complex conjugation induces a contactomorphism
\[
c\colon (S^3,\xi_{\rm st})\to (S^3,\xi_{\rm st})
\]
which reverses the orientation of $\xi_{\rm st}$. Since $c\circ
A=A^{-1}\circ c$ for every element $A$ of the group $G_{p,q}$ defined
in~\eqref{e:subgroup}, $c$ induces the contactomorphism $\overline c$.

Since $L(p,q)$ has no non--trivial real two--cohomology, $\om$ is
exact. Eliashberg~\cite[Proposition~3.1]{El2} proved that in this
situation there exist a contact form $\al\in\Om^1(L(p,q))$ defining
$\overline\xi_{\rm st}$, a constant $C>1$ and a symplectic form $\Om$
on 
\[
(0,+\infty)\x L(p,q)
\] 
which coincides with $d(t\al)$ on 
\[
[C,+\infty)\x L(p,q)
\]
and is equivalent to $\om$ near $\{1\}\x L(p,q)$.

Eliashberg's argument is the following. Since $\om|_{\overline\xi_{\rm st}}>
0$, we have
\begin{equation}\label{e:om=dal}
\om|_{\overline\xi_{\rm st}}=h d{\overline \al_0} |_{\overline\xi_{\rm st}}
\end{equation}
for some smooth function $h\colon L(p,q)\to (0,+\infty)$. Let $\al=h
{\overline \al_0}$. The restriction of $\om$ to $L(p,q)$ is of the
form $d\al+d\be$ for some $\be\in\Om^1(L(p,q))$, and it follows
from~\eqref{e:om=dal} that $d(t\al)+d\be$ is a symplectic form near
$\{1\}\x L(p,q)$. By the Symplectic Tubular Neighborhood Theorem,
there exists a self--diffeomorphism $\psi$ of a tubular neighborhood
of $\{1\}\x L(p,q)$ which is the identity restricted to $\{1\}\x
L(p,q)$ and which sends $\om$ to $d(t\al)+d\be$. Let 
\[
g\colon (0,+\infty)\to [0,1] 
\]
be a smooth function which is identically $1$ on $(0,1]$ and vanishes
on $[C,+\infty)$ for some constant $C>1$. Then, if $|g'(t)|$ is
sufficiently small the two--form 
\[
\Om=d(t\al)+d(g\be)
\] 
is a symplectic form on 
\[
(0,+\infty)\x L(p,q). 
\]
The condition on $g'(t)$ is easily fulfilled as long as $C$ is
sufficiently large, so $\Om$ is the desired symplectic form.

Now consider the diffeomorphism
\begin{eqnarray*}
F\colon & (0,+\infty)\x L(p,q) & \longrightarrow (0,+\infty )\x L(p,q)\\
& (t, x) & \longmapsto \left(t h(x) , x\right)
\end{eqnarray*}
Then, $F^*(d(t{\overline \al_0}))=d(t\al)$. By first gluing via $\psi$
near $\{1\}\x L(p,q)$ and then via $F$ near $\{C+1\}\x L(p,q)$, we
obtain $(\widetilde W,\widetilde\om)$.
\end{proof}

\begin{proof}[Proof of Proposition~\ref{p:gluing}]
The statement is an immediate consequence of Theorem~\ref{t:mcw},  
Lemma~\ref{l:compatible} and Lemma~\ref{l:end}.
\end{proof}

\begin{lem}\label{l:hamiltonian}
Let $p>q\geq 1$ be coprime integers, and suppose that 
\[
\frac p{p-q}=[b_1,\ldots,b_k],\quad b_1,\ldots, b_k\geq 2.
\]
Let $N=1+\sum_{i=1}^k (b_i-1)$. Then, there exist a symplectic form
$\om$ on 
\[
X_N:=\CP^2\# N{\overline\CP}^2 
\]
and an embedding 
\[
L(p,q)\hra X_N 
\]
which is $\om$--compatible with respect to the
$S^1$--action~\eqref{e:standardaction}, such that $X_N^-$ is a
compact, regular neighborhood of a symplectic string of type
\[
(1,1-b_1,-b_2,\ldots,-b_k)
\]
in $X_N$.
\end{lem}

\begin{proof}
The proof of the lemma is a simple adaptation 
of~\cite[Theorem~2.1 and Corollary~2.1]{McW}, therefore we only outline the
argument here, and refer to~\cite{McW} for details. The point is to
show that there exist: 
\begin{enumerate}
\item
an $S^1$--Hamiltonian symplectic form $\om$ on $X_N$, and  
\item
a Hamiltonian $H\colon X_N\to\R$ for $\om$ and a regular value $c\in\R$ of
$H$, 
\end{enumerate}
such that 
\[
X^-=\{H\geq c\}.
\]

Let $S_{-1}$ be a ruled surface over $\CP^1$ with zero section $Z_0$
and infinity section $Z_{\infty}$ satisfying $Z_0\cdot Z_0=-1$ and
$Z_\infty\cdot Z_\infty=1$. There is a standard $S^1$--action on
$S_{-1}$ which rotates the fibers and fixes $Z_0$ and $Z_\infty$, and
an $S^1$--Hamiltonian symplectic structure on $S_{-1}$ with a
Hamiltonian function $H\colon S_{-1}\to [0,1]$ such that $H(Z_0)=0$ and
$H(Z_\infty)=1$. Fix a point $x\in Z_0$. The union of $Z_0$, $Z_1$ and
a fiber of the fibration $S_{-1}\to\CP^1$ through the point $x$ is a
symplectic string of type $(-1,0,1)$. Blowing up at $x$ and taking
proper transforms yields an $S^1$--Hamiltonian symplectic structure on
$\widehat S_{-1}$ together with a symplectic string of type
$(-2,-1,-1,1)$. We can keep blowing up at the intersection points of
the proper transforms of $Z_0$ until the initial part of
the string looks like this:
\[
(-a_1,-a_2,\ldots,-a_h,-1,\ldots),
\]
where the integers $a_i\geq 2$ are chosen so that 
\[
\frac pq=[a_1,\ldots,a_h]. 
\]
Using the Riemenschneider's point rule~\cite{Ri} it is easy to check
that this can be done so that the whole string is of the following:
\[
(-a_1,-a_2,\ldots,-a_h,-1,-b_k,-b_{k-1},\ldots,-b_2,1-b_1,1).
\]
Now the resulting blowup of $S_{-1}$ is $X_N$. As explained in the
proof of~\cite[Theorem~2.1]{McW} (see also~\cite[Lemmas~2.3]{McW}
and~\cite[Lemmas~2.4]{McW}), $X_N$ carries an $S^1$--Hamiltonian
symplectic form $\om$ with Hamiltonian 
\[
H\colon X_N\to\R, 
\]
there is a regular value $c$ of $H$ such that 
\[
X^+=\{H\leq c\}
\] 
is a regular neighborhood of a symplectic string $\Ga_1$ of type
\[
(-a_1,\ldots,-a_h), 
\]
and 
\[
X^-=\{H\geq c\}
\] 
is a regular neighborhood of a symplectic string $\Ga_2$ of type
\[
(1,1-b_1,-b_2,\ldots,-b_k).
\]
Moreover, $H^{-1}(c)$ is isomorphic, as an $S^1$--manifold, to
$L(p,q)$ endowed with the $S^1$--action~\eqref{e:standardaction}. This
is because the $S^1$--action~\eqref{e:standardaction} is the action
induced on $L(p,q)$ when the lens space is viewed as the boundary of
the $S^1$--equivariant plumbing determined by the weighted graph dual
to $\Ga_1$, which coincides, by~\cite[Corollary~2.1]{McW}, with the
$S^1$--action induced on $H^{-1}(c)$.
\end{proof}

\begin{proof}[Proof of Theorem~\ref{t:string}]
By Lemma~\ref{l:hamiltonian} we can find a symplectic structure $\om$
on a rational symplectic four--manifold $X$ and an embedding
$L(p,q)\hra X$ which is $\om$--compatible with respect to the
$S^1$--action~\eqref{e:standardaction} and such that $X^-$ is a
regular neighborhood of a symplectic string
\[
\widetilde\Ga=\widetilde C_0\cup \widetilde C_1\cup\cdots \widetilde C_k 
\subset X
\]
of type $(1,1-b_1,-b_2,\ldots, -b_k)$.  Applying
Proposition~\ref{p:gluing} we can construct a symplectic manifold of
the form 
\[
(X_W=W\cup_{\del W=\del X^-} X^-,\overline\om).
\]

By~\cite{McD1}, if $(M,\om)$ is a closed symplectic four--manifold
containing an embedded symplectic two--sphere $C$ of
self--intersection $+1$ such that $M\setminus C$ is minimal, then
$(M,\om)$ is symplectomorphic to $\CP^2$ with the standard K\"ahler
form. Morever, the symplectomorphism can be chosen such that it sends
$C$ to a complex line. Since $\widetilde C_0\subset X_W$ has
self--intersection $+1$, and since non--minimal symplectic
four--manifolds can be reduced to minimal ones by blowing down
exceptional symplectic spheres, we conclude that, for some $M\geq 0$,
there is a symplectomorphism
\[
\psi\colon (X_W,\overline\om) \to (\CP^2\# M{\overline \CP}^2,\om_M)
\] 
sending $\widetilde C_0$ to a complex line $C_0\subset\CP^2$. Clearly,
\[
\Ga=C_0\cup\cdots\cup C_k=\psi(\widetilde \Ga)
\subset\CP^2\# M{\overline \CP}^2 
\]
is a symplectic string of type $(1,1-b_1,-b_2,\ldots,-b_k)$.
\end{proof}

\section{Complements of symplectic strings}
\label{s:strings}

The purpose of this section is to prove Theorem~\ref{t:main}(1), which
follows immediately from Theorem~\ref{t:string} combined with 
the following result:

\begin{theorem}\label{t:diff-uniqueness}
Let $\om_M$ be a symplectic form on $X_M=\CP^2\# M\overline\CP^2$
obtained from the standard K\"ahler form on $\CP^2$ by symplectic
blowups. Let
\[
\Ga=C_0\cup C_1\cup\cdots\cup C_k\subset X_M
\]
be a symplectic string of type 
\[
(1,1-b_1,-b_2,\ldots,-b_k),\quad k\geq 2,\quad b_1,\ldots,b_k\geq 1,
\]
and such that $C_0\subset\CP^2$ is a complex line. Let $s_i$, for
$i=1,\ldots, k$, denote one half of the cardinality of the set
\[
S_i:=\{e\in H_2(X_M;\Z)\ | \ e\cdot e = -1,\ 
e\cdot [C_i]\neq 0,\ e\cdot [C_j]=0
\ \text{for}\ j\in\{0,\ldots,k\}\setminus\{i\}\}
\]
and let $n_i:=b_i-s_i$. Then, $\n:=(n_1,\ldots,n_k)\in\bZ_{p,q}$ and
the complement of a regular neighborhood of $\Ga$ is
orientation preserving diffeomorphic to a smooth blowup of
$W_{p,q}(\n)$.
\end{theorem}

Theorem~\ref{t:diff-uniqueness} will follow from
Theorem~\ref{t:proper-transform} below.

\begin{theorem}\label{t:proper-transform}
Let $\om_M$ be a symplectic form on $X_M=\CP^2\# M\overline\CP^2$
obtained from the standard K\"ahler form on $\CP^2$ by symplectic
blowups. Let
\[
\Ga=C_0\cup C_1\cup\cdots\cup C_k\subset X_M
\]
be a symplectic string of type 
\[
(1,1-b_1,-b_2,\ldots,-b_k),\quad k\geq 2,\quad b_1,\ldots,b_k\geq 1,
\]
and such that $C_0\subset\CP^2$ is a complex line. Then, there is 
a sequence of symplectic blowdowns
\[
(\CP^2\# M\overline\CP^2,\om_M)\to (\CP^2\# (M-1)\overline\CP^2,\om_{M-1})
\to\cdots\to (\CP^2,\om_0)
\]
with $\om_0$ diffeomorphic to the standard K\"ahler form and such that
$\Ga$ descends to a symplectic string of type $(1,1)$ in $(\CP^2,\om_0)$.
\end{theorem}

Before proving Theorems~\ref{t:diff-uniqueness}
and~\ref{t:proper-transform} we need to establish some auxiliary
results.

\begin{lem}\label{l:string}
Let $\{l,f_1,\ldots,f_M\}\subseteq H_2(\CP^2\# M\overline\CP^2;\Z)$ be a
set of generators which are orthogonal with respect to the
intersection form, with $l$ being the homology class of a complex line
in $\CP^2$ and with each $f_i$ having square $-1$. Let $\al_1,\ldots,
\al_k\in H_2(\CP^2\# M\overline\CP^2;\Z)$, $k\geq 2$, be homology
classes of the form
\begin{equation*}
\al_1 = l-e^1_1-e^1_2-\cdots - e^1_{b_1},\quad
\al_i = e^i_1 - e^i_2 - \cdots - e^i_{b_i},\quad i=2,\ldots, k,
\end{equation*}
where $b_i\geq 1$ for $i=1,\ldots, k$, $e^i_j\in\{f_1,\ldots,f_M\}$
for every $i,j$ and $e^i_j\neq e^i_{j'}$ if $j\neq j'$. Suppose also
that
\[
\al_i\cdot\al_j=
\begin{cases}
1\quad\text{if $|i-j|=1$}\\
0\quad\text{if $|i-j|>1$}
\end{cases}
\]
and let $A^1=\{e^1_1,\ldots, e^1_{b_1}\}$ and $A^i=\{e^i_2,\ldots,
e^i_{b_i}\}$ for $i=2,\ldots, k$. Then, 
\begin{enumerate}
\item
for every $j=2,\ldots, k$, there is an index $i$ with $1\leq i < j$ 
such that 
\[
e^j_1\in A^i
\]
Moreover, if $e^j_1\in A^i$ with $i<j-1$, 
\[
e^h_1\in A^i\cap A^j
\]
for some index $h$ with $i<h<j$;
\item
for every $1\leq i<j\leq k$ we have
\begin{equation*}
A^i\cap A^j\subseteq\{e^2_1,\ldots,e^k_1\}.
\end{equation*}
\end{enumerate}
\end{lem}

\begin{proof}
(1) We argue by induction on $k\geq 2$. For $k=2$, (1) applies only to
$j=2$, in which case the statement clearly holds because
$\al_1\cdot\al_2=1$.

Now suppose that $k>2$ and the statement to be true for
$\al_1,\ldots,\al_{k-1}$.  Clearly, it suffices to prove the statement
for $j=k$. Since $\al_{k-1}\cdot\al_k=1$, either $e^k_1\in A^{k-1}$,
in which case we are done, or $e^{k-1}_1\in A^k$.  In the latter case,
we set $j_1=k-1$ and the induction hypothesis implies $e^{k-1}_1\in
A^{j_2}$ for some $j_2<j_1$. Thus, $A^{j_2}\cap A^k\neq\emptyset$
and, since $\al_{j_2}\cdot\al_k=0$, we either have $e_1^k\in A^{j_2}$,
in which case, setting $i=j_2$ and $h=k-1$ we are done because
$e^{k-1}_1\in A^{j_2}\cap A^k$ or $e^{j_2}_1\in A^k$. Continuing in
this fashion we obtain a maximal, strictly decreasing sequence of
indexes
\[
k=j_0>k-1=j_1>j_2>\cdots > j_r\geq 1
\]
such that $e^{j_s}_1\in A^{j_{s+1}}\cap A^k$ for $s=1,\ldots,
r-1$. Then, $e^k_1\in A^{j_r}$. In fact, if $e^k_1\not\in A^{j_r}$
then, since $\al_1\cdot\al_k=0$, we would have $j_r >1$, therefore the
sequence could be extended contradicting its maximality. This proves
(1) setting $i=j_r$ and, when $j_r < k-1$, $h=j_{r-1}$.

(2) The proof is again an induction on $k\geq 2$. For $k=2$,
$\al_1\cdot\al_2$ implies $A^1\cap A^2=\emptyset$, so suppose $k>2$
and the statement to be true for $\al_1,\ldots,\al_{k-1}$.

Observe that, since $\al_i\cdot\al_k\geq 0$ for $i<k$, 
\[
e_1^k\not\in\{e^1_1,\ldots,e_1^{k-1}\}.
\]
This implies that $e_1^k$ belongs to at most one of the
sets $A^1, A^2,\ldots, A^{k-1}$,  because $e_1^k\in A^i\cap A^j$ for
$1\leq i<j\leq k-1$ implies 
\[
A^i\cap A^j\not\subseteq\{e^2_1,\ldots,e^{k-1}_1\}, 
\]
contrary to the induction hypotheses.

Arguing by contradiction, suppose the statement to be false for
$\al_1,\ldots,\al_k$. By the induction hypothesis, for some $1\leq
i\leq k-1$ we have
\begin{equation*}
A^i\cap A^k\not\subseteq\{e^2_1,\ldots,e^k_1\}.
\end{equation*}

We claim that $i\neq k-1$. In fact, if $i=k-1$ then
$\al_{k-1}\cdot\al_k=1$ and $A^{k-1}\cap A^k\neq \emptyset$ imply 
$e^k_1\in A^{k-1}$ and $e_1^{k-1}\in
A^k$. But by (1) there is an index $i'<k-1$ such that $e^{k-1}_1\in
A^{i'}$. Thus, since $\al_{i'}\cdot\al_k=0$, $e^k_1\in A^{k-1}$, and 
by the induction hypothesis
\[
A^{i'}\cap A^{k-1}\subseteq\{e_1^2,\ldots, e_1^{k-1}\}, 
\]
we have $e^h_1\in A^k$. Now we can apply (1) again with $j=i'$ and
argue in the same way. After a finite number of similar steps we are
forced to conclude $e_1^k\in A^1$, which is incompatible with
$e^k_1\in A^{k-1}$.

We can now finish the proof assuming $i\neq k-1$. Since
$\al_i\cdot\al_k=0$, either $e^i_1\in A^k$ or $e^k_1\in A^i$. In the
latter case, by (1) we have $e^h_1\in A^i\cap A^k$ for some index $h$
with $i<h<k$. But $A^i\cap A^k\not\subseteq\{e^2_1,\ldots,e^k_1\}$
implies that $A^i\cap A^k$ has at least two elements. Therefore, in
any case we have $e^i_1\in A^k$.  Repeated applications of (1)
starting with $j=i$ and $j=k-1$ respectively, yield two maximal
strictly decreasing sequences
\begin{equation*}
i=h_0>h_1>\cdots >h_r\geq 1,\hspace{1cm} k-1=l_0>l_1>\cdots >l_s\geq 1,
\end{equation*}
such that 
\[
e_1^{h_n}\in A^{h_{n+1}}\cap A^k, \quad
n=0,\ldots, r-1,
\]
and
\[
e_1^{l_m}\in A^{l_{m+1}}\cap A^k,
\quad m=0,\ldots,s-1.
\]
Maximality implies $e_1^{h_r}\not\in A^k$if $h_r>1$,  
$e_1^{l_s}\not\in A^k$ if $l_s>1$, and $e_1^k\in
A^{h_r}\cap A^{l_s}$. Since $e_1^k$ belongs to at most one of the sets
$A^1,\ldots, A^{k-1}$, we must have $h_r=l_s$. Therefore there exist
$\overline n$ with $0\leq {\overline n}< r$ and $\overline m$ with
$0\leq {\overline m} < s$ such that $h_{\overline n}\neq l_{\overline
m}$ but $h_{{\overline n}+1}=l_{{\overline m}+1}$. This implies that
each of the distinct elements $e^{h_{\overline n}}_1$ and $e^{l_{\overline
m}}_1$ belong to $A^{h_{{\overline n}+1}}\cap A^k$. Since
$\al_{h_{{\overline n}+1}}\cdot\al_k=0$, we must have
$e_1^{h_{\overline n + 1}}\in A^k$, $e_1^k\in A^{h_{\overline n +
1}}$, and therefore $r=\overline n + 1=\overline m+1=s$, which is
incompatible with $e_1^{h_r}\not\in A^k$.
\end{proof}

\begin{prop}\label{p:hom-classes}
Let $X_M=\CP^2\# M\overline\CP^2$ be endowed with a blowup $\om_M$ of
the standard K\"ahler form on $\CP^2$, and let
\begin{equation*}
\Ga=C_0\cup C_1\cup\cdots\cup C_k \subset X_M
\end{equation*}
be a symplectic string of type 
\[
(1,1-b_1,-b_2,\ldots,-b_k),\quad b_1,\ldots,b_k\geq 1.
\]
Let $\{l,f_1,\ldots,f_M\}\subseteq H_2(\CP^2\# M\overline\CP^2;\Z)$ be
a set of generators which are orthogonal with respect to the
intersection form, with $l$ being the homology class of a complex line
in $\CP^2$ and with each $f_i$ having square $-1$. Suppose that
$[C_0]\in H_2(X_M;\Z)$ is equal to the homology class $l$. Then,
\begin{equation}\label{e:hom-classes}
[C_1] = l-e^1_1-e^1_2-\cdots - e^1_{b_1},\quad
[C_i] = e^i_1 - e^i_2 - \cdots - e^i_{b_i},\quad i=2,\ldots, k,
\end{equation}
where $e^i_j\in\{f_1,\ldots,f_M\}$ for every $i,j$ and $e^i_j\neq
e^i_{j'}$ if $j\neq j'$ .
\end{prop}

\begin{proof}
We can write
\[
[C_i]=\de_{1i}l+\sum_{j=1}^M a^i_j f_j 
\]
for some coefficients $a^i_j\in\Z$, $i=1,\ldots, k$, where $\de_{1i}$
is Kr\"onecker's delta. Since each $C_i$ is symplectic, 
\[
\langle c_1(X_M), [C_i]\rangle = 2 + C_i\cdot C_i\quad \text{for}\quad
i=1,\ldots, k,
\]
which is equivalent to:
\begin{equation}\label{e:adjunction}
\sum_{j=1}^M \left(a^i_j+(a^i_j)^2\right) = 2(1-\de_{1i}).
\end{equation}
We assume $k\geq 1$ and prove the statement by induction on
$k$. Equation~\eqref{e:adjunction} implies immediately
$a^1_j\in\{0,-1\}$ for $j=1,\ldots, M$. Therefore, the statement holds
for $k=1$. 

Now suppose that $k\geq 2$ and that the classes $[C_i]$ have the form
given in the statement for $i=1,\ldots, k-1$. By
Equation~\eqref{e:adjunction} for $i=k$, there is exactly one index
$j_0\in\{1,\ldots,M\}$ such that $a^k_{j_0}\in\{1,-2\}$, while
$a^k_j\in\{0,-1\}$ for $j\neq j_0$. 

We claim that the equality $a^k_{j_0}=-2$ leads to a contradiction. In
fact, since $C_k\cdot C_{k-1}=1$, if all the coefficients $a^k_j$ are
non--positive, then we must have $k>2$ and $a^{k-1}_1\in A^k$. By
Lemma~\ref{l:string}(1) applied for $j=k-1$, $e^{k-1}_1\in A^i$ for some
$i<k-1$. Since $C_i\cdot C_k=0$ and all the coefficients $a^k_j$ are
non--positive, we must have $i>1$ and $a^i_1\in A^k$. Now we can apply
Lemma~\ref{l:string}(1) again for $j=i$, and argue in the same
way. Clearly, after a finite number of similar steps we reach a
contradiction.
\end{proof}

\begin{lem}\label{l:sphere}
Under the assumptions of Theorem~\ref{t:proper-transform}, let $J$ be
an almost complex structure tamed by $\om_M$ and such that $\Ga$ is
$J$--holomorphic. Then, there exists an embedded $J$--holomorphic
sphere $\Si\subset X_M$ such that $[\Si]\cdot [C_0]=0$ and
$[\Si]\cdot[\Si]=-1$. Moreover, there exists such a $\Si$ disjoint
from $\Ga$ if and only if there exists a symplectic sphere $S\subset
X_M$ of square $-1$ such that $[S]\cdot [C_i]=0$ for $i=0,\ldots,k$.
\end{lem}

\begin{proof}
Since the symplectic 4--manifold $X_M$ is obtained by blowing up
$\CP^2$, there exists a symplectic sphere $S\subset X_M$ of square
$-1$ orthogonal to $l=[C_0]$. By~\cite[Lemma~3.1]{McD1}, the homology
class $[S]$ is either represented by an embedded $J$--holomorphic
sphere $\Si$ or by a $J$--holomorphic~\emph{cusp--curve}
\[
S_1\cup\cdots\cup S_n, 
\]
i.e.~a union of (not necessarily embedded) $J$--holomorphic
spheres. In the first case, the first part of the lemma is proved. In
the second case notice that, since $[S]\cdot [C_0]=0$, by positivity
of intersections~\cite{McD2} we have
\[
[S_i]\cdot [C_0]=0\quad\text{for}\quad i=1,\ldots, n.
\]
Therefore 
\[
[S_i]\cdot [S_i]\leq -1\quad\text{for}\quad i=1,\ldots, n.
\]
Moreover,
\[
1 = \chi(X) + S\cdot S = \langle c_1(X_M), [S]\rangle = 
\sum_{i=1}^n \langle c_1(X_M), [S_i]\rangle = 
\sum_{i=1}^n (\chi(S_i) + [S_i]\cdot [S_i]),
\]
which implies $[S_j]\cdot [S_j]=-1$ for at least one index
$j\in\{1,\ldots, n\}$. By the adjunction formula~\cite{McD1},
$\Si:=S_j$ is embedded. Hence, the first part of the lemma is proved.
If $[S]$ is orthogonal to all the classes $[C_i]$, then by positivity
of intersections so is $[\Si]$, and therefore $\Si$ must be disjoint
from $\Ga$. This proves the second part of the lemma.
\end{proof}

\begin{proof}[Proof of Theorem~\ref{t:proper-transform}]
Let $J$ be an almost complex structure tamed by $\om_M$ and such that
$\Ga$ is $J$--holomorphic. If there exists an embedded symplectic
sphere $S\subset X_M$ such that $[S]\cdot [C_i]=0$, $i=0,\ldots,k$,
then by Lemma~\ref{l:sphere} there is an embedded $J$--holomorphic
sphere $\Si\subset X_M\setminus\Ga$ with self--intersection
$\Si\cdot\Si=-1$. Therefore, we may blow down $\Si$ and reapply the
same argument to $X_{M-1}$. After a finite number of similar steps we
get $\Ga\subset X_{M'}$, with the property that for every symplectic
sphere $S\subset X_M$ the class $[S]$ intersects non--trivially at
least one of the classes $[C_i]$. Applying Lemma~\ref{l:sphere} again,
we know that there exists an embedded $J$--holomorphic sphere
$\Si\subset X_{M'}$ such that $[\Si]\cdot [C_0]=0$ and
$\Si\cap\Ga\neq\emptyset$. By Proposition~\ref{p:hom-classes}, the
homology classes $[C_i]$ have the form~\eqref{e:hom-classes}. Since
$[\Si]\cdot [C_0]=0$, $[\Si]$ must coincide with one of the classes
$e^i_j$.

Now we argue by induction on $k\geq 2$. Clearly, either for some
$s\geq 2$ we have $b_s=1$ and $[\Si]=[C_s]$, or $[\Si]\cdot
C_i\in\{0,1\}$ for every $i=1,\ldots, k$.

If $b_s=1$ and $[\Si]=[C_s]$, by positivity of intersections we must
have $\Si=C_s$. In this case we can blow down $\Si$, and $\Ga$
descends to a symplectic string $\Ga_1\subset X_{M'-1}$ of length
strictly less than the length of $\Ga$.  If $k=2$ then $s=2$,
$b_1=b_2=1$ and $\Ga_1=C_0\cup C'_1$ is a symplectic string of type
$(1,1)$ with $[C'_1]=l$. Since the complement of $C_0$ is minimal, the
conclusion follows immediately from the results of~\cite{McD1} as in the
proof of Theorem~\ref{t:string}. If $k>2$, since the intersection form of
$X_{M'-1}$ restricted to the orthogonal complement of the class $l$ is
negative definite, $\Ga_1$ satisfies the hypothesis of the theorem. By
induction, the statement holds for $\Ga_1$ and therefore for $\Ga$.

If $[\Si]\cdot C_i\in\{0,1\}$ for every $i=1,\ldots,k$, we must have
$\Si\cdot C_i=1$ for at least one index $i>0$. By
Lemma~\ref{l:string}(2), in this case there is exactly one such index,
so $\Si$ must intersect $\Ga\setminus C_0$ transversely at one smooth
point. As before, if we blow down $\Si$ then $\Ga$ descends to a
symplectic string $\Ga_1\subset X_{M'-1}$ satisfying the hypothesis of
the theorem. Now we can go through the same process starting from the
beginning, thus showing that $\Ga_1$ descends to a symplectic string
$\Ga_2\subset X_{M'-2}$. After a finite number of similar steps we
arrive at a symplectic string of length strictly less than the length
of $\Ga$, and the induction argument works as before.
\end{proof}

\begin{proof}[Proof of Theorem~\ref{t:diff-uniqueness}]
By Theorem~\ref{t:proper-transform}, there is 
a sequence of symplectic blowdowns
\begin{equation}\label{e:blowdowns}
(\CP^2\# M\overline\CP^2,\om_M)\to (\CP^2\# (M-1)\overline\CP^2,\om_{M-1})
\to\cdots\to (\CP^2,\om_0)
\end{equation}
with $\om_0$ diffeomorphic to the standard K\"ahler form, and such
that $\Ga$ descends to a symplectic string of type $(1,1)$ in
$(\CP^2,\om_0)$. Since we want to determine the complement of a
regular neighborhood of $\Ga$ up to diffeomorphisms, we may assume
that $\om_0$ is the standard K\"ahler form on $\CP^2$.

Let $J_0$ be an almost complex structure in $\CP^2$ tamed by $\om_0$
and such that $l_0$ and $l'_0$ are $J_0$--holomorphic. The results
of~\cite{Gr} show that if $\{J_t\}_{t\in [0,1]}$ is a generic path of
tamed almost complex structures connecting $J_0$ to the standard
complex structure $J_1$, then there are families $\{l_t\}_{t\in
[0,1]}$ and $\{l'_t\}_{t\in [0,1]}$ of smooth, embedded and distinct
$J_t$--holomorphic spheres connecting $l_0$, respectively $l'_0$, to
standard complex lines $l_1$ and $l'_1$. This shows that the string
$l_0\cup l'_0$ is symplectically isotopic to a pair $l_1\cup l'_1$ of
distinct standard complex lines.

Thus, in order to determine the diffeomorphism type of the complement
of a regular neighborhood of $\Ga$, one may replace $\Ga$ with the
proper transform of a pair of distinct complex lines in $\CP^2$ under
the sequence of complex blowups corresponding to~\eqref{e:blowdowns}.
In fact, by analyzing the construction of Sequence~\eqref{e:blowdowns}
in the proof of Theorem~\ref{t:proper-transform} and using Kirby
calculus and Lemma~\ref{l:0-sequence}, it is easy to see that there is an
orientation preserving diffeomorphism
\[
W\cong W_{p,q}(\n)\# r\overline\CP^2,
\]
with $\n\in\bZ_{p,q}$ as in the statement and $r=M-\sum_{i=1}^k
(b_i-n_i)$.
\end{proof}

\section{Stein structures on $W_{p,q}(\n)$}
\label{s:stein}

In this section we construct Stein structures on the smooth
four--manifolds with boundary $W_{p,q}(\n)$ defined in
Section~\ref{s:intro}. The proof is based on Legendrian
surgery~\cite{El1, Go}. 

A knot $K$ in a contact three--manifold $(Y,\xi)$ is
called~\emph{Legendrian} if $K$ is everywhere tangent to the
distribution $\xi$. The contact structure induces a framing of $K$,
usually called the~\emph{contact framing} and denoted by
$\tb(K)$. Assume that $\xi$ is oriented. Given an oriented Legendrian
knot $K$ in $(Y,\xi)$ and a non--zero section $v$ of $\xi$ along $K$,
the~\emph{rotation number} $\rot_v(K)$ is the winding number of the
oriented tangent vector to $K$ with respect to $v$ in the oriented
plane $\xi$.

Let $(n_1,\ldots, n_k)$ be an admissible sequence of positive integers
such that 
\[
[n_1,\ldots,n_k]=0. 
\]
We fix once and for all a sequence of strict blowdowns
\begin{equation}\label{e:sequence}
(n_1,\ldots, n_k)\to\cdots\to (0)
\end{equation}
 as in Lemma~\ref{l:0-sequence}. We can realize this sequence of
operations geometrically, by viewing each step as a blowdown in the
sense of the Kirby calculus on framed links, starting from the framed
link $L=\cup_{i=1}^k L_i$ of Figure~\ref{f:L} and ending with the
zero--framed unknot. Such a sequence corresponds to an
orientation preserving diffeomorphism
\begin{equation}\label{e:varphi}
\varphi\colon N(\n)\to S^1\x S^2, 
\end{equation}
where $N(\n)=N(n_1,\ldots, n_k)$ is the closed oriented
three--manifold obtained by surgery along the framed link $L\subset
S^3$. Let $\nu(L)\subset S^3$ be a small tubular neighborhood of
$L$. Then, the complement $S^3\setminus \nu(L)$ can be identified with
a subset $\CC$ of $N(\n)$, i.e.~the complement of the surgered solid
tori. Every link in $\CC\subset S^3$ is therefore endowed with a
canonical framing. We shall use this canonical framing to identify any
framing with a $k$--uple of integers.

The smooth manifolds $W_{p,q}(\n)$ are obtained by first attaching a
one--handle to the four--ball, and then attaching two--handles to the
boundary of the resulting $S^1\x D^3$. The standard tight contact
structure $\ze_0$ on $S^1\x S^2=\del (S^1\x D^3)$ is obtained from the
standard contact structure $\xi_{\rm st}$ on $S^3$ by removing two
smooth balls and gluing the resulting boundaries in a suitable
way. Moreover, each Legendrian link in $(S^1\x S^2, \ze_0)$ is
contact isotopic to a Legendrian link in \emph{standard form} in the
sense of~\cite{Go}. In particular, there is a nowhere vanishing
section of $\ze_0$, denoted by $\frac{\del}{\del x}$ in~\cite{Go},
along any Legendrian link in standard form. The section
$\frac{\del}{\del x}$ extends as a nowhere vanishing section $v$ of
$\ze_0$ to all of $S^1\x S^2$. The contact structure $\ze_0$ has a
natural orientation coming from the natural complex orientation
on $\xi_{\rm st}$. Therefore, the rotation numbers with respect to the
section $v$ are well--defined.

Using Eliashberg's Legendrian surgery construction~\cite{El1,Go}, we will
now prove that there are Stein structures on each $W_{p,q}(\n)$ by
showing that the attaching circles of the two--handles are isotopic to
Legendrian knots 
\[
K_i\subset (S^1\x S^2,\ze_0)
\]
and the two--handles are attached with framings $\tb(K_i)-1$.

Let $\widetilde\ze_0=\varphi^*(\ze_0)$ be the tight contact structure
on $N(\n)$ obtained by pulling back $\ze_0$ via the
diffeomorphism~\eqref{e:varphi}. Observe that $\widetilde\ze_0$ does
not depend on the choice of $\varphi$ because $S^1\x S^2$ carries only
one tight contact structure up to isotopy. Let $\tau$ be the
pull--back of the nowhere zero section $v$ by the
diffeomorphism~\eqref{e:varphi}.

Recall from Section~\ref{s:intro} that, given two coprime integers
\[
p>q\geq 1,\quad\text{with}\quad 
\frac{p}{p-q}=[b_1,\ldots,b_k], 
\]
we defined the set
\[
\bZ_{p,q}=\{(n_1,\ldots,n_k)\in\Z^k\ |\ [n_1,\ldots,n_k]=0, 
\  0\leq n_i\leq b_i,\ i=1,\ldots, k\}.
\]

\begin{figure}[ht]
  \setlength{\unitlength}{1mm}
  \begin{center}
    \begin{picture}(100,30)
      \put(7,0){\large $\L_1$}
      \put(23,0){\large $\L_2$}
      \put(68,0){\large $\L_{k-1}$}
      \put(84,0){\large $\L_k$}
      \put(10,30){\large $n_1$}
      \put(26,30){\large $n_2$}
      \put(68,30){\large $n_{k-1}$}
      \put(85,30){\large $n_k$}
      \put(0,3){\includegraphics[width=10cm]{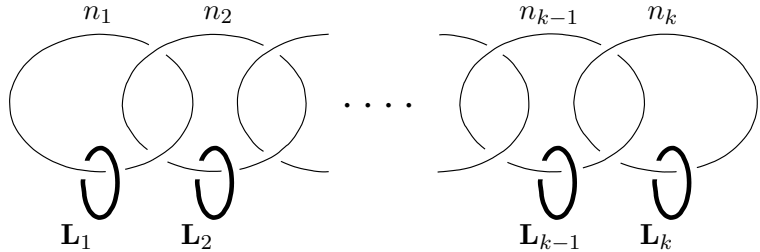}}
    \end{picture}    
  \end{center}
\caption{The link ${\mathbf L}$} 
\label{f:link} 
\end{figure}

\begin{theorem}\label{t:legendrian}
Let $p>q\geq 1$ be coprime integers with
$\frac{p}{p-q}=[b_1,\ldots,b_k]$, and let
$\n=(n_1,\ldots,n_k)\in\bZ_{p,q}$.  Let 
\[
{\mathbf L}=\cup_{i=1}^k {\mathbf L}_i\subset
\CC=S^3\setminus\nu(L) 
\]
be the ``thick'' link drawn in Figure~\ref{f:link}.
Then, there exists a Legendrian link 
\[
\LL=\bigcup_{i=1}^k \LL_i\subset (\CC,\widetilde\ze_0|_{\CC}) 
\]
with the following properties:
\begin{enumerate}
\item[(a)]
$\LL$ is smoothly isotopic to ${\mathbf L}$ inside $\CC$.
\item[(b)] 
$\LL$ has contact framing equal to the canonical framing induced
by the inclusion $\CC\subset S^3$.
\item[(c)] Define $\rot_\tau(\LL_i)$ to be $0$ for $i<1$ and
$i>k$. Then, $\LL$ admits an orientation such that
$\rot_\tau(\LL_1)=0$, and
\[
\rot_\tau(\LL_{i-1}) + \rot_\tau(\LL_{i+1}) - n_i\rot_\tau(\LL_i) = 
\begin{cases}
\phantom{-}1\quad \text{if $i=1$},\\
\phantom{-}0\quad\text{if $1<i<k$}\\
-1\quad\text{if $i=k$}
\end{cases}
\]
 \end{enumerate}
\end{theorem}
\begin{proof}
We argue by induction on $k\geq 1$. When $k=1$ we have
$\frac{p}{p-q}=[b_1]$ and $[n_1]=0$, therefore $b_1=p$, $q=p-1$,
$n_1=0$ and $\varphi$ is the identity. The formulas for the
Thurston--Bennequin and the rotation number of a Legendrian link in
standard form~\cite{Go} show that the Legendrian link $\LL$ of
Figure~\ref{f:L(p,p-1)} satisfies (a), (b) and (c).
\begin{figure}[ht]
\begin{center}
\includegraphics[height=1in]{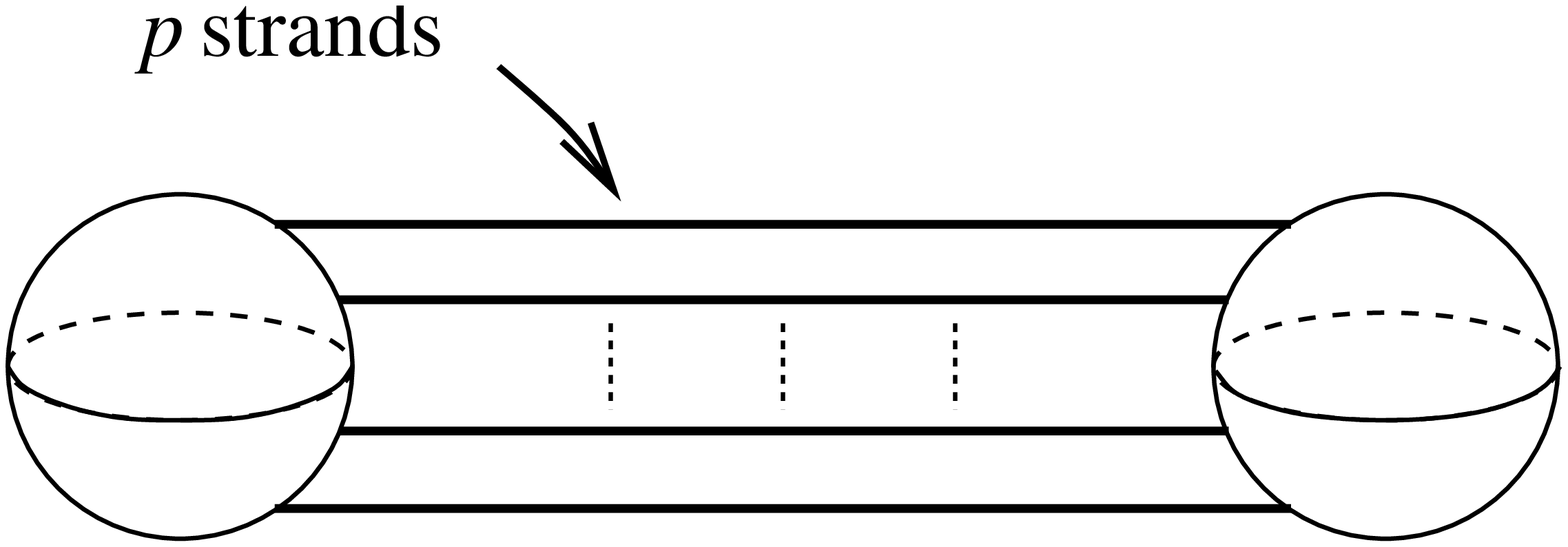} 
\end{center}
\caption{The link $\LL$ when $k=1$}
\label{f:L(p,p-1)} 
\end{figure}
Now assume $k>1$, and suppose that the statement holds for every
four--manifold of the form $W_{p',q'}(\m)$, where
$\m=(m_1,\ldots,m_{k-1})$.  Let
\begin{equation*}
\n=(n_1,\ldots, n_{s-1}, 1, n_{s+1},\ldots, n_k)\to 
\n'=(n_1,\ldots, n_{s-1}-1, n_{s+1}-1,\ldots, n_k)
\end{equation*} 
be the first element of the sequence of strict
blowdowns~\eqref{e:sequence}. By definition of a strict blowdown,
$s>1$.  Redraw Figure~\ref{f:link} as Figure~\ref{f:braid}, where the
thin horizontal arcs in the picture are the strands of the braid
$\be=\si_1^2\si_2^2\cdots\si_{k-1}^2\in B_k$. Observe that the closure
of $\be$ is isotopic to the thin link $L$ of Figure~\ref{f:link}. We
orient ${\mathbf L}$ as shown in Figure~\ref{f:braid}.
\begin{figure}[ht]
\setlength{\unitlength}{1mm}
\begin{center}
    \begin{picture}(120,60)
      \put(10,56){\large $\L_1$}
      \put(15.5,43){\large $\L_2$}
      \put(10,37){\large $\L_{s-1}$}
      \put(17,31.2){\large $\L_s$}
      \put(23,19){\large $\L_{s+1}$}
      \put(9,15){\large $\L_{k-1}$}
      \put(15,0){\large $\L_k$}
      \put(110,55){\large $n_1$}
      \put(110,50){\large $n_2$}
      \put(110,36){\large $n_{s-1}$}
      \put(110,31){\large $1$}
      \put(110,27){\large $n_{s+1}$}
      \put(110,13){\large $n_{k-1}$}
      \put(110,8){\large $n_k$}
      \put(0,3){\includegraphics[width=12cm]{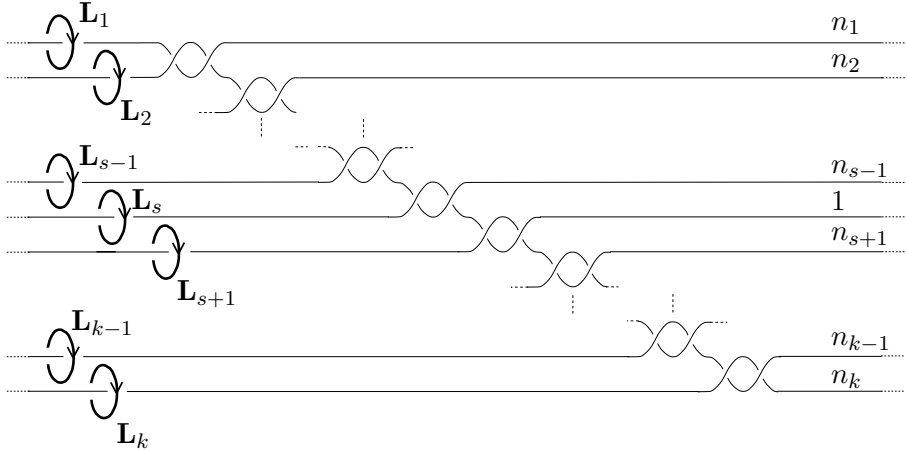}}
    \end{picture}    
 \end{center}
\caption{The link ${\mathbf L}$ redrawn and oriented}
\label{f:braid} 
\end{figure}

We now prove the statement assuming $s\neq k$, which is the harder
case. At the end we will briefly say how to deal with the easier case
$s=k$, omitting some obvious details. 

Blowing down the $s$--th component of the thin link $L$ yields
Figure~\ref{f:blowndown}.
\begin{figure}[ht]
\setlength{\unitlength}{1mm}
  \begin{center}
    \begin{picture}(120,60)
      \put(9,53){\large $\L'_1$}
      \put(15,38){\large $\L'_2$}
      \put(9,34){\large $\L'_{s-1}$}
      \put(69,33){\large $\L''_s$}
      \put(17,19){\large $\L'_{s+1}$}
      \put(9,15){\large $\L'_{k-1}$}
      \put(16,0){\large $\L'_k$}
      \put(108,49.5){\large $n_1$}
      \put(108,45){\large $n_2$}
      \put(108,31){\large $n_{s-1}-1$}
      \put(108,26.7){\large $n_{s+1}-1$}
      \put(108,13){\large $n_{k-1}$}
      \put(108,8.2){\large $n_k$}
      \put(0,3){\includegraphics[width=12cm]{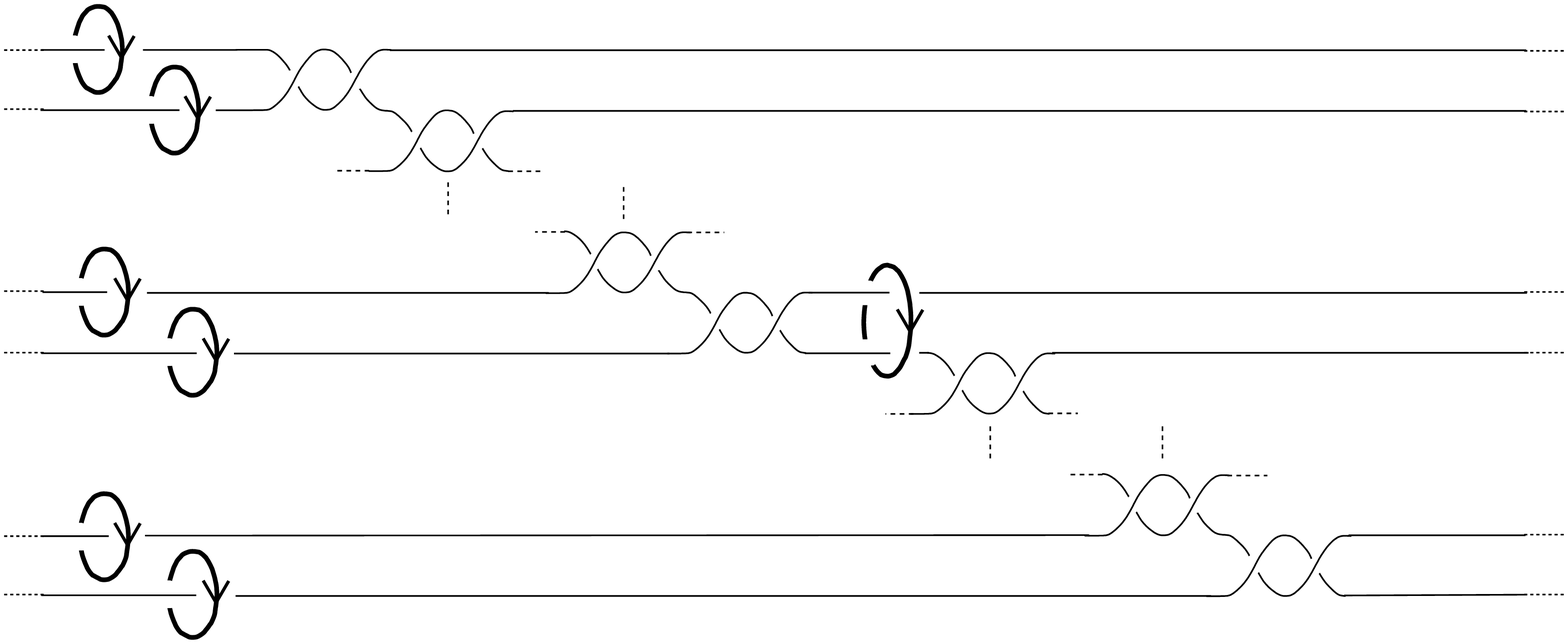}}
    \end{picture}    
  \end{center}
\caption{The link ${\mathbf L}$ blown down} 
\label{f:blowndown} 
\end{figure}
The resulting thick link $\widetilde{\mathbf L}$ is the disjoint union
of ${\mathbf L}'$ and ${\mathbf L}''_s$, where ${\mathbf
L}'=\bigcup_{i\neq s}{\mathbf L}'_i$ is the collection of thick
unknots, each of which is linked to a single strand of the thin braid,
on the left--hand side of the picture, and ${\mathbf L}''_s$ is the
only thick unknot which links two strands. We can view the link
$\widetilde{\mathbf L}$ as sitting inside the oriented three--manifold
$N(\n')$ obtained by surgery in $S^3$ along the thin framed link of
Figure~\ref{f:blowndown}, which we call $L'$. The diffeomorphism
$\varphi$ of~\eqref{e:varphi} is the composition of two
diffeomorphisms:
\[
N(\n)\stackrel{\psi}{\lra} 
N(\n')\stackrel{\varphi'}{\lra}
S^1\x S^2.
\]
Here $\psi$ is determined by the first element in the
sequence~\eqref{e:sequence}, and $\varphi'$ by the remaning
elements. Let $\widetilde\ze_0'$ be the pull--back to $N(\n')$ of
$\ze_0$ via $\varphi'$. Observe that the oriented knot ${\mathbf
L}''_s$ is isotopic to an oriented band connected sum of the oriented
knots ${\mathbf L}'_{s-1}$ and ${\mathbf L}'_{s+1}$ as shown in
Figure~\ref{f:band}.
\begin{figure}[ht]
\setlength{\unitlength}{1mm}
  \begin{center}
    \begin{picture}(110,45)
      \put(24,16){\large $B$}
      \put(19,34){\large $\L'_{s-1}$}
      \put(76,34){\large $\L''_s$}
      \put(19,0){\large $\L'_{s+1}$}
      \put(95,30){\large $n_{s-1}-1$}
      \put(95,9){\large $n_{s+1}-1$}
      \put(0,2){\includegraphics[width=11cm]{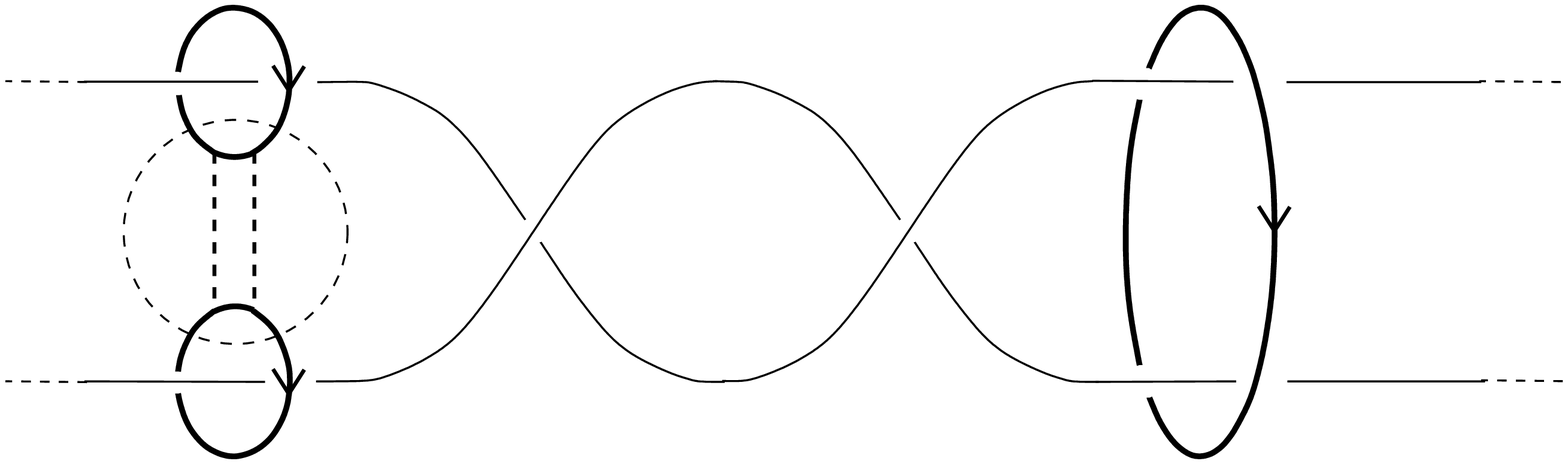}}
    \end{picture}    
  \end{center}
\caption{The band connected sum} 
\label{f:band} 
\end{figure}
By the induction hypothesis, the oriented link ${\mathbf
L}'=\cup_{i\neq s}{\mathbf L}'_i$ is isotopic in the complement of
$L'$ to a link which is Legendrian with respect to
$\widetilde{\ze'_0}$ and satisfies conditions (b) and (c) of the
statement. Therefore, without loss of generality we may assume
${\mathbf L}'$ to be a Legendrian link $\LL'=\bigcup_{i\neq s}\LL'_i$,
with each of its components having contact framing equal to the
canonical framing, and the rotation numbers with respect to the
section $v'=(\varphi')^{-1}_*(\tau)$ satisfying the stated relations
(c).

Since any two sufficiently small Legendrian arcs in a contact
three--manifold are contact isotopic (see e.g.~\cite{EF}), without
loss of generality we may assume that (i) there is a contactomorphism
$g$ between the dotted three--ball $B$ in Figure~\ref{f:band} and a
ball centered at the origin of $\R^3$ endowed with the standard
contact structure $\{dz+xdy=0\}$ and (ii) the intersection of $B$ with
$\LL'_{s-1}$ and $\LL'_{s+1}$ is sent by the contactomorphism onto two
horizontal arcs sitting in the $yz$--plane, one above and the other
below the $xy$--plane, as in the left--hand side of
Figure~\ref{f:leg1}. The right--hand side of Figure~\ref{f:leg1}
describes, in the language of front projections (cf.~\cite{Go}), how
to perform a Legendrian band connected sum of $\LL'_{s-1}$ and
$\LL'_{s+1}$. The result is a Legendrian knot $\LL''_s$ smoothly
isotopic to ${\mathbf L}''_s$ in $\CC$.
\begin{figure}[ht]
\begin{center}
\includegraphics[height=1.7in]{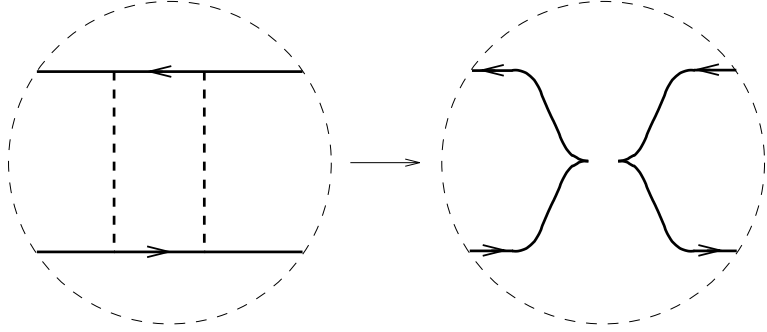} 
\end{center}
\caption{The Legendrian band connnected sum} 
\label{f:leg1} 
\end{figure}
Using the standard formula computing the Thurston--Bennequin number of
a Legendrian knot from its front projections, it is easy to check that
$\LL''_s$ has contact framing equal to its canonical framing minus
one. On the other hand, the diffeomorphism $\psi^{-1}$ sends the
canonical framing of $\LL''_s$ to the canonical framing plus one,
while the canonical framings of ${\mathcal L}'_i$, for $i\neq s$, are
sent by $\psi^{-1}$ to the canonical framings. This implies that each
component of ${\mathcal L}=\psi^{-1}({\mathcal L}'\cup\LL'_s)$ has
canonical framing equal to the contact framing with respect to
$\widetilde\ze_0=\psi^{-1}(\widetilde\ze'_0)$, concluding the proof of
(a) and (b) when $s\neq k$.

To prove (c) when $s\neq k$, let $v'=(\varphi'_*)^{-1}(v)$. Induction
applied to the oriented Legendrian link $\LL'$ inside $N(\n')$ gives
the relations:
\begin{equation}\label{e:relation1}
\rot_{v'}(\LL'_{i-1}) + 
\rot_{v'}(\LL'_{i+1}) - 
n_i \rot_{v'}(\LL'_i) =
\begin{cases}
\phantom{-}1 \quad i=1\\
\phantom{-}0 \quad 1<i<k,\quad i\neq s, s\pm 1\\
-1 \quad i=k
\end{cases}
\end{equation}
and
\begin{equation}\label{e:relation1a}
\rot_{v'}(\LL'_{s-2}) + 
\rot_{v'}(\LL'_{s+1}) - 
(n_{s-1}-1) \rot_{v'}(\LL'_{s-1}) =
\begin{cases}
1 \quad s=2\\
0 \quad 2<s<k
\end{cases}
\end{equation}
\begin{equation}\label{e:relation1b}
\rot_{v'}(\LL'_{s-1}) + 
\rot_{v'}(\LL'_{s+2}) - 
(n_{s+1}-1) \rot_{v'}(\LL'_{s+1}) =
\begin{cases}
\phantom{-}0 \quad 1<s<k-1\\
-1 \quad s=k-1
\end{cases}
\end{equation}
On the other hand, by Figure~\ref{f:leg1} and the formula for the
rotation number of an oriented Legendrian knot in terms of its front
projections~\cite{Go}, we have
\begin{equation}\label{e:relation2}
\rot_{v'}(\LL'_{s-1}) + 
\rot_{v'}(\LL'_{s+1}) - 
\rot_{v'}(\LL''_s) = 0
\end{equation}
Now set $\LL_i=\psi^{-1}(\LL'_i)$ for $i\neq s$, and
$\LL_s=\psi^{-1}(\LL''_s)$. The relations~\eqref{e:relation1},
\eqref{e:relation1a}, \eqref{e:relation1b} and~\eqref{e:relation2}
provide, when pulled--back via $\psi$, the stated relations for
$\LL=\cup_i\LL_i$. This concludes the proof when $s\neq k$.

When $s=k$ the argument is similar but simpler. The main difference is
that Figure~\ref{f:leg1} should be replaced with
Figure~\ref{f:leg2}. The rest of the argument is essentially the same,
so we omit the details.
\begin{figure}[ht]
\begin{center}
\includegraphics[height=1.7in]{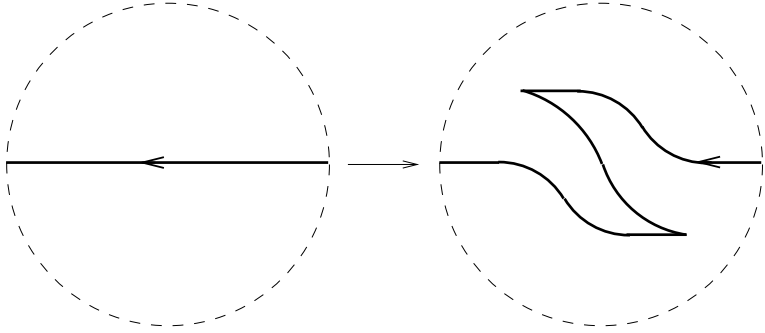} 
\end{center}
\caption{The case $s=k$} 
\label{f:leg2} 
\end{figure}
\end{proof}

\begin{cor}\label{c:stein}
Let $p>q\geq 1$ be coprime integers with
$\frac{p}{p-q}=[b_1,\ldots,b_k]$, and let
$\n=(n_1,\ldots,n_k)\in\bZ_{p,q}$. Fix a diffeomorphism $\varphi$ as
in~\eqref{e:varphi}, and let $\LL=\cup_i\LL_i$ be an oriented
Legendrian link as in Theorem~\ref{t:legendrian}. Let
$\widetilde\LL\subset (N(\n),\widetilde\ze_0)$ be a Legendrian link
consisting of $b_i-n_i$ distinct Legendrian push--offs of each
component $\LL_i$ of $\LL$, for $i=1,\ldots,k$. Then, $W_{p,q}(\n)$
carries a Stein structure obtained by performing Legendrian surgery
along the Legendrian link
\[
\varphi(\widetilde\LL)\subset (S^1\x S^2,\ze_0).
\]
\end{cor}
 
\begin{proof}
The statement follows immediately from Theorem~\ref{t:legendrian}
using Legendrian surgery~\cite{El1,Go}.
\end{proof}

\section{Contact structures on $\del W_{p,q}(\n)$}
\label{s:contact}

In this section we prove Theorem~\ref{t:main}(2). The main ingredient
of the proof will be the following result.

\begin{theorem}\label{t:part1}
Let $p>q\geq 1$ be coprime integers and $\n\in\bZ_{p,q}$.  Let $J$ be
a Stein structure on $W_{p,q}(\n)$ constructed as in
Corollary~\ref{c:stein}, and let $\om$ be a symplectic form on
$W_{p,q}(\n)$ compatible with $J$. Then, $(W_{p,q}(\n),\om)$ is a
symplectic filling of $(L(p,q), {\overline\xi_{\rm st}})$.
\end{theorem}

Given a Stein structure $J$ on $W_{p,q}(\n)$, the distribution 
\[
\xi =T\del W_{p,q}(\n)\cap J T\del W_{p,q}(\n)
\]
of complex lines tangent to the boundary is a contact
structure. Therefore, if $\om$ is a symplectic form on $W_{p,q}(\n)$
compatible with $J$, then $(W_{p,q}(\n),\om)$ is a symplectic filling
of $(\del W_{p,q}(\n),\xi)$. We shall establish Theorem~\ref{t:part1}
by proving that if $J$ is constructed as in Corollary~\ref{c:stein},
then the contact three--manifold $(\del W_{p,q}(\n),\xi)$ is
isomorphic to $(L(p,q),\overline\xi_{\rm st})$.

The contact structure $\xi$ is tight~\cite{El3}. Therefore, in view of
the classification of the tight contact structures on
$L(p,q)$~\cite{Gi, Ho}, to show that $(\del W_{p,q}(\n),\xi)$ and
$(L(p,q),\overline\xi_{\rm st})$ are isomorphic it suffices prove
that, after a suitable identification $\del W_{p,q}(\n)=L(p,q)$, $\xi$
and $\overline\xi_{\rm st}$ induce the same Spin$^c$ structure on
$L(p,q)$, i.e. are homotopic as two--plane fields in the complement of
a point. In order to do this we shall use Gompf's
$\Ga$--invariant.

Let $M$ be a closed, oriented three--manifold and denote by $\SS(M)$
the set of spin structures on $M$. Given an oriented two--plane field
$\xi$ on a $M$, Gompf~\cite{Go} defines a map
\[
\Ga(\xi,\cdot)\colon\SS(M)\to H_1(M;\Z)
\]
which depends only on the homotopy class $[\xi]$. Moreover, reversing
the orientation of $M$ reverses the sign of $\Ga(\xi,s)$ for fixed
$s$, and fixing $q\in M$, $\Ga(\cdot,s)$ classifies oriented
two--plane fields on $M\setminus\{q\}$ up to homotopy.

As shown in~\cite[Theorem~4.12]{Go}, it is possible to compute $\Ga$
when $\xi$ is the distribution of complex lines tangent to the
boundary of an almost complex four--manifold. We will need a slight
generalization of that result.

Let $X^*$ be a smooth four--manifold obtained by attaching
two--handles to $B^4$ along a framed link $\La\subset S^3$. Let
$\La_1\subset \La$ be a sublink, and let $X^*_1\subset X^*$ be the
submanifold obtained by attaching the two--handles corresponding to
$\La_1$. We may think of $X^*$ as obtained by attaching two--handles
along $\La_2=\La\setminus \La_1$ to the boundary of $X_1^*$. Suppose
that the boundary of $X_1^*$ is a connected sum of $S^1\x S^2$'s, so
that $\del X_1^*=\del X_1$, where $X_1$ is obtained by attaching
one--handles to $B^4$. Then, the smooth four--manifold
$X=X_1\cup\left(X^*\setminus X_1^*\right)$ is well--defined because
every self--diffeomorphism of $\del X_1$ extends to $X_1$. Suppose
that $X$ carries an almost complex structure $J$, and denote by $\xi$
the distribution of complex lines tangent to $\del X$. Fix a complex
trivialization $\tau$ of $TX$ over $X_1$. Restricting to $\del X_1$,
$\tau$ determines a Spin structure $s_0$ on $\del X_1$, and since the
boundaries of $X_1$ and $X_1^*$ are identified, there is a canonically
associated characteristic sublink $\La_0\subset \La_1$ representing
the Poincar\'e dual to the second Stiefel--Whitney class of $TX^*_1$
relative to $s_0$~\cite{GS}. Similarly, to any Spin structure $s$ on
$\del X=\del X^*$ one can canonically associate a characteristic
sublink $\La(s)\subset\La$. Choose an orientation for $\La$, let
$K_1,\ldots, K_n$ be its components and $\al_1,\ldots, \al_n\in
H_2(X^*;\Z)$ be the corresponding two--homology classes.  The
following statement can be proved by an almost word--for--word
repetition of the proof of~\cite[Theorem~4.12]{Go}.
\begin{theorem}[\cite{Go}]\label{t:gompf-quoted}
The Poincar\'e dual to $\Ga(\xi,s)$ is equal to the restriction of the
class $\rho\in H^2(X^*;\Z)$ determined by the evaluations:
\begin{equation*}
\langle \rho, \alpha_i\rangle = \frac{r_i + \lk(K_i, \La_0 + \La(s))}2,\quad 
i=1,\ldots, n
\end{equation*}
where $r_i$ is equal to zero if $K_i\subset \La_1$, and to the integer
obstruction to extending the trivialization $\tau$ over the
corresponding two--handle if $K_i\subset \La_2$.
\qed\end{theorem}

\begin{rem}\label{r:aftergompfquoted}
(a) Theorem~\ref{t:gompf-quoted} reduces to~\cite[Theorem~4.12]{Go}
 when $X_1$ is endowed with a standard Stein structure, $X$ is
 obtained by Legendrian surgery along a Legendrian link $\La_2$ in
 ``standard form'', and the restriction of $\tau$ to $\del X_1$ is
 induced by the vector field ``$\frac{\del}{\del x}$''
 (cf.~\cite{Go}). In this case each component $K_i$ of the oriented
 link $\La_2$ has a well--defined rotation number -- the relative
 winding number of an oriented tangent vector to $K_i$ with respect to
 $\frac{\del}{\del x}$ -- which coincides with the number $r_i$ as
 defined above.

(b) As observed in~\cite{Go}, after the statement of Theorem~4.12,
replacing $\La_0+\La(s)$ in the formula by any smooth one--cycle
carried by $\La$ and agreeing with $\La_0+\La(s)$ modulo $2$ does not
change the restriction of $\rho$. Thus, after such a replacement
Theorem~\ref{t:gompf-quoted} still holds.

(c) Theorem~\ref{t:gompf-quoted} applies to another particular case,
i.e.~when $\La_1=\emptyset$ and $X^*_1=X_1=B^4$. In this case
$\La_0=\emptyset$ and for every $K_i\subset \La_2$ we have
$r_i=\langle c_1(J), \al_i\rangle$.
\end{rem}

We are going to apply Theorem~\ref{t:gompf-quoted} when the link $\La$
is the framed link $L\cup\L\subset S^3$ of Figure~\ref{f:W}, with
$\La_1=L$ and $\La_2=\L$, respectively, the ``thin'' and the ``thick''
link. We choose an orientation for $\La$ as follows. Represent
$L\cup\L$ as in Figure~\ref{f:braid}. Orient $\L$ as shown in the
picture. Orient $L$ that each of its components has linking number
$+1$ with one of the corresponding thick meridian shown in
Figure~\ref{f:braid}. 

In this case $X_1$ is $S^1\x D^3$ and $X$ is $W_{p,q}(\n)$.  Let $s$
be a Spin structure on $\del W_{p,q}(\n)$ determined by a
characteristic sublink $\La(s)\subset \La$. Observe that $s$ is
uniquely specified by a sequence $(s_1,\ldots,s_k)$, where
$s_i\in\{0,1\}$ is equal to $1$ if the $i$--th component of $L$
belongs to $\La(s)$ and to $0$ otherwise. In fact, since $\La(s)$ is
characteristic, it is easy to check that the components of $\L$
belonging to $\La(s)$ are determined by the components of $L$ belong
to $\La(s)$.

Denote by $\mu_i\in H_1(\del W_{p,q}(\n);\Z)$, $i=1,\ldots, k$,
homology classes corresponding to positively oriented meridians of the
components of $L$. The classes $\{\mu_i\}$ generate the first homology
of $\del W_{p,q}(\n)$ and satisfy the relations
\begin{equation}\label{e:mu-relations}
b_i \mu_i = \mu_{i-1} + \mu_{i+1}, \quad i=1,\ldots, k,
\end{equation}
where $\mu_i$ is to be interpreted as the zero class for $i<0$ and
$i>k$.

\begin{prop}\label{p:gamma-invariant}
Let $p>q\geq 1$ be coprime integers, and let $\n\in\bZ_{p,q}$.  Let
$J$ be a Stein structure on $W_{p,q}(\n)$ constructed as in
Corollary~\ref{c:stein}, and let $\xi$ be the contact structure
induced on the boundary. Let $s$ be the Spin structure on $\del
W_{p,q}(\n)$ determined by the sequence $(s_1,\ldots, s_k)$, and set
$s_i=0$ for $i<1$ and $i>k$. Then,
\begin{equation*}
\PD\Ga(\xi,s) = \sum_{i=1}^k \frac{s_{i-1}+s_{i+1} + b_i(1-s_i)}2
\mu_i -\sum_{i=2}^k \mu_i
\end{equation*}
\end{prop}

\begin{proof}
Let $W^*$ be the four--manifold with boundary obtained by attaching
two--handles to $B^4$ according to the framed link $\La=L\cup{\mathbf
L}$ of Figure~\ref{f:W}. According to Theorem~\ref{t:legendrian},
${\mathbf L}$ is isotopic in the complement of $L$ to a link
${\mathcal L}$, Legendrian with respect to the contact structure
$\widetilde\ze_0$. 

Let $s_0$ be the unique spin structure on $S^1\x S^2$ which extends
over a two--handle attached along $S^1\x\{p\}$ if and only if the
corresponding framing is even. Fix a diffeomorphism $\varphi$ as
in~\eqref{e:varphi}, and let $\tilde s_0$ be the pull--back of $s_0$
under $\varphi$. Let $L_0\subset L$ be the characteristic sublink
corresponding to $\tilde s_0$.

Let $\al_1,\ldots,\al_k\in H_2(W^*;\Z)$ be the homology classes
determined by the components $L_1,\ldots, L_k$ of $L$, and let
$\{\be^{j_i}_i\ |\ i=1,\ldots, k,\ j_i=1,\ldots, b_i-n_i\}\subset
H_2(W^*;\Z)$ be the classes determined by the components ${\mathcal
L}_i^{j_i}$ of ${\mathcal L}$. Let $L(s)\subset\La$ be the
characteristic sublink corresponding to the Spin structure $s$
determined by the sequence $(s_1,\ldots,s_k)$. Then, by
Theorem~\ref{t:gompf-quoted}, $\PD\Ga(\xi,s)$ is equal to the
restriction to $\del W^*=\del W_{p,q}(\n)$ of the class $\rho\in
H^2(W^*;\Z)$ determined by the values:
\begin{align*}
\langle\rho,\al_i\rangle &= \frac12 \lk(L_i,L_0+L(s)),\\
\langle\rho,\be^{j_i}_i\rangle &= \frac12 (r_i + 
\lk({\mathcal L}^{j_i}_i,L_0+L(s))).
\end{align*}
Notice that
\[
 \frac12 (r_i + \lk({\mathcal L}^{j_i}_i,L_0+L(s)))=
 \frac12 (r_i + \lk({\mathcal L}_i,L_0+L(s)))
\]
for every $j_i=1,\ldots,b_i-n_i$, where ${\mathcal L}_i$ is the
$i$--th component of the link $\mathcal L$ given in
Theorem~\ref{t:legendrian}.  Setting 
\[
C_i=\frac12 \lk(L_i,L_0+L(s)),\qquad
D_i=\frac12 (r_i + \lk({\mathcal L}_i,L_0+L(s))), 
\]
we have:
\begin{equation}\label{e:gamma-invariant}
\PD\Ga(\xi,s)= \sum_{i=1}^k \left(C_i + D_i (b_i-n_i)\right) \mu_i.
\end{equation}
Observe that when $b_i=n_i$ there is no component $\LL^{j_i}_i$, but
the number $D_i$ is still well--defined. Hence,
Equation~\eqref{e:gamma-invariant} holds in any case. If we extend the
definition of $C_i$ and $D_i$ by setting $C_i=D_i=0$ for $i<0$ and
$i>k$, in view of the relations~\eqref{e:mu-relations} we have:
\begin{equation}\label{e:gamma-invariant2}
\begin{split}
\PD\Ga(\xi,s) 
=& \sum_{i=1}^k C_i\mu_i + \sum_{i=1}^k D_i (\mu_{i-1}+\mu_{i+1}-n_i\mu_i) \\
=&\sum_{i=1}^k (C_i + D_{i-1} + D_{i+1} - n_i D_i)\mu_i.
\end{split}
\end{equation}
The following identities are easy to check:
\[
\lk(\LL_{i-1}, L_j) + \lk(\LL_{i+1}, L_j) -
n_i\lk(\LL_i, L_j) + \lk(L_i, L_j) = 0,\quad i,j=1,\ldots,k
\]\[
\lk(\LL_i, L(s)) = 2s_i - 1, \quad i=1,\ldots,k
\]\[
\lk(L_i, L(s)) = b_i-n_i+s_i(2n_i-b_i)-s_{i-1}-s_{i+1},\quad
i=1\ldots, k.
\]
Here a linking number is to interpreted as zero if there is an index
less than $1$ or bigger than $k$. In view of
Theorem~\ref{t:legendrian}(c), a simple calculation using the identities
above gives:
\[
2(C_i+D_{i-1}+D_{i+1}-n_i D_i)=
\begin{cases}
s_2+s_{0}+b_1(1-s_1), \quad i=1\\
s_{i+1}+s_{i-1}+b_i(1-s_i)-2, \quad i=2,\ldots,k
\end{cases}
\]
The statement follows immediately substituting these values
in~\eqref{e:gamma-invariant2}.
\end{proof}

\begin{proof}[Proof of Theorem~\ref{t:part1}]
Let $\xi$ be the contact structure on $\del W_{p,q}(\n)$ given by the
tangent complex lines. As explained at the beginning of the section,
in order to prove the rest of the statement it suffices to show that,
after the choice of a suitable identification $\del
W_{p,q}(\n)=L(p,q)$, for any Spin structure $s$ on $L(p,q)$ we have:
\[
\Ga(\xi, s)=\Ga({\overline\xi_{\rm st}}, s).
\]
Recall that $(L(p,q),{\overline\xi_{\rm st}})$ is the link of a cyclic
quotient singularity. Let $R_{p,q}$ a regular neighborhood of the
exceptional divisor inside the canonical resolution of such a
singularity. Then $R_{p,q}$ is diffeomorphic to a plumbing of type
$(-a_1,\ldots,-a_h)$, with
\begin{equation}\label{e:c1values}
\langle c_1(R_{p,q}), x_i\rangle = 2 - a_i, \quad i=1,\ldots, h, 
\end{equation}
where $x_i\in H_2(R_{p,q};\Z)$ is the obvious generator represented by
an embedded rational curve $C_i$ with self--intersection $-a_i$ (see
e.g.~\cite{BPV}). Equations~\eqref{e:c1values} follow from the
adjunction formula.

If we start with an immersed curve in $\CP^2$ which is the union
$l_1\cup l_2$ of two distinct complex lines, we can successively blow
up the curve and its proper transforms at points not of $l_1$, until
we obtain a string $C$ of curves of type $(1,1-b_1,-b_2,\ldots,-b_k)$
inside $\CP^2\# N\overline\CP^2$. As the proof of
Lemma~\ref{l:hamiltonian} shows, there is a natural
orientation preserving diffeomorphism $\varphi$ between the
complement $Z$ of a regular neighborhood $\nu(C)$ of $C$ and
$R_{p,q}$. Moreover, since the natural generators of $H_2(Z;\Z)$ are
complex curves, by the adjunction formula $\varphi$ must preserve the
first Chern classes. Since $R_{p,q}$ is simply connected, this implies
that $\varphi$ preserves the complex structures up to homotopy. It
follows that the distribution of complex lines tangent to $\del Z$ is
homotopic to ${\overline\xi_{\rm st}}$.

The oriented three--manifold $\del\nu(C)$ has a surgery presentation
given by a chain of $k+1$ unknots $U_0,\ldots, U_k$ framed,
respectively, $1,1-b_1,\ldots,-b_k$. A Spin structure on $\nu(C)$ is
encoded by a characteristic sublink of $\cup_{i=0}^k U_i$
(see~\cite[Section~5.7]{GS}), which can be identified with a
$(k+1)$--tuple $(t_0,t_1,\ldots,t_k)$, $t_i\in\{0,1\}$, by requiring
that $U_j$ belongs to the sublink if and only if $t_j=1$.

Blowing down $U_0$ gives an identification
\[
\nu(C)=-L(p,q)=L(p,p-q). 
\]
By looking at the effect of the Kirby move on the Spin
structure~\cite[pp.~190--191]{GS}, one sees that $t_0=1-t_1$.

We shall compute
\[
\Ga_{L(p,q)}(\overline\xi_{\rm
st},s)=-\Ga_{L(p,p-q)}(\overline\xi_{\rm st},s)
\]
applying Theorem~\ref{t:gompf-quoted} to the almost complex
four--manifold $\nu(C)$.  

Define
\begin{equation*}
b'_i=
\begin{cases}
1 & i=0\\
-b_1+1 & i=1\\
-b_i & i=2,\ldots,k
\end{cases}
\end{equation*} 
By Remark~\ref{r:aftergompfquoted}(b) and (c), we can apply
Theorem~\ref{t:gompf-quoted} with $L_0=\emptyset$ and the sign of
$L(s)$ reversed. There are natural generators 
\[
\nu_0,\ldots,\nu_k\in H_1(\del\nu(C);\Z) 
\]
corresponding to the meridians of $U_0,\ldots,
U_k$, such that
\begin{equation}\label{e:finalformula}
\begin{split}
\PD\Ga_{L(p,p-q)}(\overline\xi_{\rm st},s) = 
\sum_{i=0}^k \frac{2 + b'_i (1-s_i) - t_{i-1} - t_{i+1}}2 \nu_i =\\ 
\frac{2+(1-t_0)-t_1}2 \nu_0 + \frac{2-(b_1-1)(1-t_1) -t_0 - t_2}2 \nu_1\\ 
+\sum_{i=2}^k \frac{2-b_i(1-t_i) -t_{i-1} -t_{i+1}}2 \nu_i
\end{split}
\end{equation}

Blowing down all the $(-1)$--unknots in Figure~\ref{f:W} identifies
the boundary of $W_{p,q}(\n)$ with $L(p,q)=-\del\nu(C)$. Under this
identification a Spin structure $(s_1,\ldots,s_k)$ on $\del
W_{p,q}(\n)$ corresponds to 
\[
(t_0,t_1,\ldots,t_k)=(1-s_1,s_1,\ldots,s_k), 
\]
and each generator $\mu_i$ is sent to $\nu_i$, $i=1,\ldots,k$.

Since $\nu_0=-\nu_1$ and $t_0=1-t_1$, Equation~\eqref{e:finalformula}
and Proposition~\ref{p:gamma-invariant} give
\[
\PD\Ga_{L(p,q)}(\overline\xi_{\rm st},s)=\PD\Ga_{L(p,q)}(\xi,s). 
\]
This concludes the proof.
\end{proof}

The following is a restatement of Theorem~\ref{t:main}(2).

\begin{cor}\label{c:(2)}
For every $\n\in\bZ_{p,q}$, $W_{p,q}(\n)$ carries a symplectic form
$\om$ such that $\left(W_{p,q}(\n),\om\right)$ is a symplectic filling
of $\left(L(p,q), {\overline\xi_{\rm st}}\right)$. Moreover, the
homology group $H_2(W_{p,q}(\n);\Z)$ contains no classes of 
self--intersection $-1$.
\end{cor}

\begin{proof}
The first part of the statement follows from Theorem~\ref{t:part1},
because a Stein four--manifold with boundary is well--known to carry a
symplectic form compatible with the complex structure. Therefore, by 
Theorems~\ref{t:part1} and~\ref{t:string} there exists a symplectic
string $\Ga$ of type $(1,1-b_1,-b_2,\ldots, -b_k)$ inside a rational
symplectic four--manifold $X_M$ such that $W_{p,q}(\n)$ is
orientation preserving diffeomorphic to the complement
$X_M\setminus\nu(\Ga)$. If $H_2(W_{p,q}(\n);\Z)$ contained a class of
square $-1$, then by Lemma~\ref{l:sphere} $W_{p,q}(\n)$ would contain
a smooth $(-1)$--sphere. But by Corollary~\ref{c:stein} $W_{p,q}(\n)$
carries a Stein structure, and a Stein four--manifold does not contain
embedded $(-1)$--spheres~\cite[Proposition~2.2]{LM2}. This proves the
second part of the statement.
\end{proof}

\section{The proof of Theorem~\ref{t:main}}
\label{s:mainproof}

Let $p>q\geq 1$ be coprime integers, with
$\frac{p}{p-q}=[b_1,\ldots,b_k]$, and let $R_{p,p-q}$ be a regular
neighborhood of the exceptional divisor in the canonical resolution of
the cyclic quotient singularity of type $(p,p-q)$.  $R_{p,p-q}$ is
diffeomorphic to a plumbing of $2$--disk bundles over $2$--spheres of
type $(-b_1,\ldots,-b_k)$, and there is a natural identification
$\del R_{p,p-q}=L(p,p-q)$.  By~\cite{Bo}, given an
orientation preserving diffeomorphism
\[
f\colon L(p,p-q)\to L(p,p-q),
\]
the isotopy class of $f$ is uniquely determined by the induced
homomorphism 
\[
f_*\colon H_1(L(p,p-q);\Z)\lra H_1(L(p,p-q);\Z). 
\]
Moreover, $f_*$ can only be multiplication by $1$, $-1$, $q$ or $-q$,
and if $f_*$ is multiplication by $\pm q$ then $q^2\equiv 1\bmod
p$. Let $x_1,\ldots, x_k\in H_2(R_{p,p-q};\Z)$ be the natural
generators represented by smooth rational curves with
self--intersections $-b_1,\ldots,-b_k$.

\begin{lem}\label{l:extension}
Let $p>q\geq 1$ be coprime integers. Let $f\colon\del R_{p,p-q}\to\del
R_{p,p-q}$ be an orientation preserving diffeomorphism such that the
induced homomorphism
\[
f_*\colon H_1(\del R_{p,p-q};\Z)\to H_1(\del R_{p,p-q};\Z)
\]
is multiplication by $q$. Then, $f$ is the restriction of a
diffeomorphism
\[
F\colon R_{p,p-q}\to R_{p,p-q}
\]
such that
\[
F_*(x_i)=x_{k+1-i},\quad i=1,\ldots, k.
\]
\end{lem}

\begin{proof}
It suffices to prove that if $q^2\equiv 1\bmod p$, then the
four--manifold $R_{p,p-q}$ admits a self--diffeomorphism which
preserves the orientation and induces multiplication by $q$ on
$H_1(L(p,p-q);\Z)$.

The condition $q^2\equiv 1\bmod p$ is equivalent to 
\[
(b_1,b_2,\ldots, b_k)=(b_k,b_{k-1},\ldots,b_1).
\]
Therefore, when $q^2=1\bmod p$ the pair $(S^3,L)$ has an involution
which can be visualized (after an isotopy of $L$) as a $\pi$--rotation
around an axis perpendicular to the plane of the picture. The
resulting self--diffeomorphism of $R_{p,p-q}$ induces multiplication
by $q$ on $H_1(L(p,p-q);\Z)$. To see this, observe that
$H_1(L(p,p-q);\Z)$ is generated by the classes $\mu_1,\ldots,\mu_k$ of
oriented meridians of $L_1,\ldots,L_k$ satisfying
Relations~\eqref{e:mu-relations} with $-b_i$ instead of $b_i$, for
$i=1,\ldots, k$. The relations imply that $\mu_1$ is a generator,
and $\mu_k=q\mu_1$, and the diffeomorphism sends $\mu_1$ to $\mu_k$.
\end{proof}

\begin{lem}\label{l:W-involution}
Let $p>q\geq 1$ be coprime integers and $\n\in\bZ_{p,q}$. For every
integer $r\geq 0$, the smooth four--manifold $\widehat
W_{p,q}(\n)=W_{p,q}(\n)\#r\overline\CP^2$ admits an
orientation preserving self--diffeomorphism which induces
multiplication by $-1$ on $H_1(\del W_{p,q}(\n);\Z)$.
\end{lem}

\begin{proof}
Rotation by $\pi$ around a horizontal axis going through each
component of the thin link in Figure~\ref{f:W} induces an
orientation preserving diffeomorphism
\[
\psi\colon W_{p,q}(\n)\lra W_{p,q}(\n).
\]
Blowing down all the $(-1)$--framed unknots gives in identification
\[
\del W_{p,q}(\n) = L(p,q),
\]
where $L(p,q)$ is viewed as the boundary of $R_{p,p-q}$ with reversed
orientation. It follows that $H_1(\del W_{p,q}(\n);\Z)$ is generated
by meridians $\mu_1,\ldots,\mu_k$ of the components of the thin link
in Figure~\ref{f:W} satifying Relations~\eqref{e:mu-relations}.
Since $\mu_1$ is a generator and $\psi$ sends $\mu_1$ to $-\mu_1$, it
follows that $\psi$ induces multiplication by $-1$ on $H_1(\del
W_{p,q}(\n);\Z)$. The same conclusion applies to $\widehat
W_{p,q}(\n)$, because we can add $r$ disjoint, unlinked
$(-1)$--framed unknots $U_1,\ldots, U_r$ to Figure~\ref{f:W} and
repeat the above argument. The fact that the link $\bigcup_i U_i$ is
invariant up to isotopy under the $\pi$--rotation implies that $\psi$ 
extends to $\widehat W_{p,q}(\n)$. 
\end{proof}

\begin{proof}[Proof of Theorem~\ref{t:main}]
Theorem~\ref{t:diff-uniqueness} and Corollary~\ref{c:(2)} give,
respectively, (1) and (2). Therefore we only need to prove (3).

If $(p',s)=(p,r)$ and $(q',\n')=(q,\n)$, then clearly 
\begin{equation}\label{e:Wdiffeo}
W_{p,q}(\n)\# r\overline\CP^2 \cong W_{p',q'}(\n')\# s\overline\CP^2.
\end{equation}
If $(p',s)=(p,r)$ and $(q',\n')=({\overline q},{\overline\n})$, we
have 
\[
\frac{p}{p-q}=
[b_1,\ldots,b_k],\quad
\frac{p}{p-q'}=[b'_1,\ldots,b'_k]
=[b_k,\ldots,b_1]. 
\]
It follows from the definitions that
\[
W_{p,{\overline q}}({\overline\n})=W_{p,q}(\n),
\]
therefore~\eqref{e:Wdiffeo} still holds. 

Conversely, suppose that~\eqref{e:Wdiffeo} holds. Then, by Part (2) we
have $r=s$. Also,~\eqref{e:Wdiffeo} implies that $L(p,q)\cong
L(p',q')$, and therefore $p'=p$ and either $q'=q$ or $q'=\overline q$.

Let us first suppose that $q'=q$. By Theorem~\ref{t:string} there
exists a rational symplectic four--manifold
\[
X_M\cong\CP^2\#M\overline\CP^2
\]
and a symplectic string
\[
\Ga=C_0\cup C_1\cup\cdots\cup C_k\subset X_M
\]
of type 
\[
(1,1-b_1,-b_2,\ldots,-b_k) 
\]
such that 
\[
W_{p,q}(\n)\# r\overline\CP^2 
\]
is orientation preserving diffeomorphic to the complement
$X_M\setminus\nu(\Ga)$ of a regular neighborhood of $\Ga$ in
$X_M$. The same holds for 
\[
W_{p,q}(\n')\# r\overline\CP^2
\]
with respect to a symplectic string
\[
\Ga'=C'_0\cup C'_1\cup\cdots\cup C'_k\subset X_{M'}
\]
of the same type. Since~\eqref{e:Wdiffeo} holds, we
must have $M=M'$ and an induced diffeomorphism
\[
\psi\colon X_M\setminus\nu(\Ga)\to X_M\setminus\nu(\Ga').
\]
Up to composing $\psi$ with the automorphism of
\[
X_M\setminus\nu(\Ga')\cong W_{p,q}(\n)\# r\overline\CP^2
\]
of Lemma~\ref{l:W-involution}, we may assume that the induced
homomorphism $\psi_*^\del$ on the first homology of the boundary is
either the identity or multiplication by $q$. Observe that
\[
\nu(\Ga)\cong\nu(\Ga')\cong\CP^2\# R_{p,p-q}. 
\]
By Lemma~\ref{l:extension}, the diffeomorphism induced on the boundary
by $\psi$ can be extended to a diffeomorphism $F\colon R_{p,p-q}\to
R_{p,p-q}$. Then, it is easy to see that F can be extended to a
diffeomorphism
\[
\widehat F\colon \CP^2\# R_{p,p-q}\to \CP^2\# R_{p,p-q}
\]
such that, under the natural isomorphism  
\[
H_2(\CP^2\# R_{p,p-q};\Z)= H_2(\CP^2;\Z)\oplus H_2(R_{p,p-q};\Z), 
\]
\[
\widehat F_*=
\begin{pmatrix}
{\rm Id} & 0\\
0 & F_*
\end{pmatrix}.
\]
It follows that $\psi$ is the restriction of a diffeomorphism
\[
\widehat\psi\colon X_M\to X_M
\]
such that if $\psi_*^\del$ is the identity then 
\begin{equation}\label{e:hom-conditions0}
\widehat\psi_* [C_i] = [C'_i],\quad i=0,\ldots, k,
\end{equation}
if $\psi_*^\del$ is multiplication by $q$, then $q^2\equiv 1\bmod p$ and
\begin{equation}\label{e:hom-conditions}
\begin{split}
\widehat\psi_* [C_0] &=[C'_0],\quad 
\widehat\psi_* [C_1] =[C'_k]+[C'_0],\\
\widehat\psi_* [C_i] &=[C'_{k-i+1}],\ i=2,\ldots,k-1,\quad  
\widehat\psi_* [C_k] =[C'_1]-[C'_0].
\end{split}
\end{equation}
By Theorem~\ref{t:diff-uniqueness}, $\n=(n_1,\ldots,n_k)$ is uniquely
determined by the homology classes $[C_1],\ldots, [C_k]$, and the same
holds for $\n'$ and the classes $[C'_i]$. Therefore, if
Equations~\eqref{e:hom-conditions0} hold then
\[
b_i-n_i=b_i-n'_i,\quad i=1,\ldots,k,
\]
hence $\n'=\n$. If $q^2\equiv 1\bmod p$ and
Equations~\eqref{e:hom-conditions} hold, then
\[
b_i-n_i=b'_{k+1-i}-n'_{k+1-i}=b_i-n'_{k+1-i},
\]
therefore $\n'=\overline\n$. This shows that if $q'=q$, then either
$(q',\n')=(q,\n)$ or $q^2\equiv 1\bmod p$ and
$(q',\n')=(q,\overline\n)$, i.e.~$(q',\n')=({\overline
q},\overline\n)$, which is what we needed to prove. 

If $q'={\overline q}$ then, since 
\[
W_{p,{\overline q}}({\n'})\# r\overline\CP^2 =
W_{p,q}({\overline\n'})\# r\overline\CP^2,
\]
by~\eqref{e:Wdiffeo} and the case $q'=q$ just proved,
$(q,{\overline\n'})$ is equal to either $(q,\n)$ or (when $q^2\equiv
1\bmod p$) $({\overline q},\overline\n)$, i.e.~$(q',\n')$ is equal to
either $({\overline q},\overline\n)$ or $(q,\n)$.
\end{proof}

\bibliographystyle{amsplain}

\end{document}